\documentclass{amsart}
\usepackage{showkeys}
\usepackage{amsmath} 
\usepackage{amssymb}
\usepackage{mathrsfs}

\newtheorem{theorem}{Theorem}[section] 
\newtheorem{claim}[theorem]{Claim}
\newtheorem{lemma}[theorem]{Lemma} 
 
\newtheorem{observation}[theorem]{Observation} 
 
\newtheorem{conclusion}[theorem]{Conclusion}

\theoremstyle{definition}
\newtheorem{definition}[theorem]{Definition}
\newtheorem{dc}[theorem]{Definition/Claim}

\newtheorem{convention}[theorem]{Convention}

\newtheorem{explanation}[theorem]{Explanation}

\newtheorem{thesis}[theorem]{Thesis}
\newtheorem{hypothesis}[theorem]{Hypothesis}

\theoremstyle{remark}
\newtheorem{question}[theorem]{Question}
\newtheorem{remark}[theorem]{Remark}
\newtheorem{notation}[theorem]{Notation}

\newcommand{\rest}{{\restriction}}
 
\newcommand{\Dom}{{\rm Dom}} 
\newcommand{\Ord}{{\rm Ord}} 
\newcommand{\Min}{{\rm Min}} 
\newcommand{\Fil}{{\rm Fil}} 
\newcommand{\FIL}{{\rm FIL}} 
\newcommand{\ZF}{{\rm ZF}} 
\newcommand{\ZFC}{{\rm ZFC}} 
\newcommand{\AC}{{\rm AC}} 
\newcommand{\DC}{{\rm DC}} 
\newcommand{\hrtg}{{\rm hrtg}} 
\newcommand{\dual}{{\rm dual}} 
\newcommand{\otp}{{\rm otp}} 
\newcommand{\wlor}{{\rm wlor}} 

\newcommand{\Rang}{{\rm Rang}}
\newcommand{\cd}{{\rm cd}}
\newcommand{\rk}{{\rm rk}}
\newcommand{\qu}{{\rm qu}}

\newcommand{\wilog}{{\rm without loss of generality}}

\newcommand{\then}{{\underline{then}}}
\newcommand{\when}{{\underline{when}}}
\newcommand{\Then}{{\underline{Then}}}

\newcommand{\Iff}{{\underline{iff}}}
\newcommand{\mn}{{\medskip\noindent}}
\newcommand{\sn}{{\smallskip\noindent}}

\newcommand{\bfv}{{\bold v}}
\newcommand{\bfw}{{\bold w}}

\newcommand{\cA}{{\mathscr A}}

\newcommand{\cB}{{\mathscr B}}

\newcommand{\gb}{{\mathfrak b}}
\newcommand{\gh}{{\mathfrak h}}
\newcommand{\ga}{{\mathfrak a}}
\newcommand{\gc}{{\mathfrak c}}

\newcommand{\bbR}{{\mathbb R}}

\newcommand{\cD}{{\mathscr D}}
\newcommand{\gd}{{\mathfrak d\/}} 

\newcommand{\cH}{{\mathscr H}}

\newcommand{\cF}{{\mathscr F}}
\newcommand{\cG}{{\mathscr G}}
\newcommand{\cI}{{\mathscr I}}

\newcommand{\cP}{{\mathscr P}}
\newcommand{\varp}{{\varepsilon}}

\newcommand{\cT}{{\mathscr T}}
 
\newcommand{\gy}{{\mathfrak y}}
\newcommand{\gx}{{\mathfrak x}}
\newcommand{\gz}{{\mathfrak z}}  
\newcommand{\cU}{{\mathscr U}}

\newcommand{\cf}{{\rm cf}}

\newcount\skewfactor
\def\mathunderaccent#1#2 {\let\theaccent#1\skewfactor#2
\mathpalette\putaccentunder}
\def\putaccentunder#1#2{\oalign{$#1#2$\crcr\hidewidth
\vbox to.2ex{\hbox{$#1\skew\skewfactor\theaccent{}$}\vss}\hidewidth}}

\newenvironment{PROOF}[2][\proofname.]
   {\begin{proof}[#1]\renewcommand{\qedsymbol}{\textsquare$_{\rm #2}$}}
   {\end{proof}}

\begin{document}

\title {pcf without choice}
\author {Saharon Shelah}
\address{Einstein Institute of Mathematics\\
Edmond J. Safra Campus, Givat Ram\\
The Hebrew University of Jerusalem\\
Jerusalem, 91904, Israel\\
 and \\
 Department of Mathematics\\
 Hill Center - Busch Campus \\ 
 Rutgers, The State University of New Jersey \\
 110 Frelinghuysen Road \\
 Piscataway, NJ 08854-8019 USA}
\email{shelah@math.huji.ac.il}
\urladdr{http://shelah.logic.at}
\thanks{This research was supported by the United States-Israel Binational
Science Foundation. 
I would like to thank Alice Leonhardt for the beautiful
typing. References like \cite[Th0.2=Ly5]{Sh:950} means the label of
Th.0.2 is y5.  The reader should note that the version in my website
is usually more udpated than the one in the mathematical archive. 
 First Typed - 2004/Jan/20. Paper 835}
 


\subjclass{Primary 03E17; Secondary: 03E05, 03E50}

\keywords {set theory, weak axiom of choice, pcf}

\date{January 18, 2019}

\begin{abstract}
We mainly investigate models of set theory with restricted
choice, e.g., ZF + DC +
the family of countable subsets of $\lambda$ is well ordered for
every $\lambda$ (really local version for a given $\lambda$). We think
that in this frame much of pcf theory, (and combinatorial set theory 
in general) can be generalized.  We prove here, in particular, 
that there is a class of regular cardinals, every large 
enough successor of singular is not measurable 
and we can prove cardinal inequalities.
\end{abstract}

\maketitle
\numberwithin{equation}{section}
\setcounter{section}{-1}
\newpage

\centerline {Anotated Content} 
\bigskip

\noindent
\S0 \quad Introduction, pg.\pageref{0}

\S(0A) \quad Background, aims and results, pg.\pageref{0A}

\S(0B) \quad Preliminaries, pg.\pageref{0B}
\mn
\begin{enumerate}
\item[${{}}$]   [Include quoting \cite{Sh:497},
(\ref{0.1},\ref{0.2}); \hrtg$(Y),\wlor(Y)$, (\ref{0.3});
defining rk$_D(f)$, (\ref{0.A} + \ref{0.B}) on $J[f,D]$,
(\ref{0.B1}, \ref{0.A}, \ref{0.D}); ${\cH}_{< \kappa,\gamma}(Y)$,
(\ref{0.E} and observation \ref{x.5}); and on closure operations (\ref{0.F}).]
\end{enumerate}
\bigskip

\noindent
\S1 \quad Representing ${}^\kappa \lambda$, pg.\pageref{1}
\mn
\begin{enumerate}
\item[${{}}$]  [We define Fil$^\ell_\kappa$ and prove a
representation theorem for ${}^\kappa \lambda$.  Essentially under ``reasonable
choice" the set ${}^\kappa \lambda$ is the union of few well ordered sets,
i.e., their number depends on $\kappa$ only".  We end by a claim on
$\Pi{\ga}$.]
\end{enumerate}
\bigskip

\noindent
\S2 \quad No decreasing sequence of subalgebras, pg.\pageref{2}
\mn
\begin{enumerate}
\item[${{}}$]  [As suggested in the title we weaken the axioms.  We
deal with ${}^\kappa \lambda$ with $\lambda^+$ not 
measurable, existence of ladder $\bar C$ witnessing cofinality
and prove that many $\lambda^+$ are regular (\ref{s.2}).]
\end{enumerate}
\bigskip

\noindent
\S3 \quad Concluding remarks, pg.\pageref{3}
\mn
\begin{enumerate}
\item[${{}}$]  [We prove that if $\mu > \kappa = \cf(\mu) > \aleph_0$,
  then from a well ordering of $\cP(\cP(\kappa)) \cup {}^{\kappa >}
  \mu$ we can define a well ordering of ${}^\kappa \mu$, see
  \ref{c2}.  If e.g. $\mu$ is strong limit singular of uncountable
  cofinality, using a well order of $\cH(\mu)$ we can define a well
  ordering of $\cP(\mu)$ hence of $H(\mu^+)$, see \ref{c4}.  Lastly,
  we give sufficient conditions (in $\ZF + \DC$) for singular $\mu$,
  that $\mu^+$ is regular, see \ref{c11}.  Actually if $\mu =
\mu^{\aleph_0} + 2^{2^\kappa},\kappa = \kappa^{\aleph_0}$ and
  $X \subseteq \mu$ codes $\cP(\cP(\kappa))$ and ${}^\omega \mu$, then
  using $X$ as a parameter we can define a well ordering of ${}^\kappa
  \mu$, see \ref{c14}.]
\end{enumerate}
\newpage

\section {Introduction}  \label{0}
\bigskip

\noindent
\subsection{Background, aims and results} \label{0A} \
\bigskip

The thesis of \cite{Sh:497} was that pcf theory without full choice
exists.  Two theorems supporting this thesis were proved.  The 
first (\cite[4.6,pg.117]{Sh:497}, we shall not mention $\ZF$) is:
\begin{theorem}
\label{0.1}  
{\rm [DC]}  If ${\cH}(\mu)$ is well ordered,
$\mu$ strong limit singular of uncountable cofinality \then \,
$\mu^+$ is regular not measurable (and $2^\mu$ is an $\aleph$,
i.e. $\cP(\mu)$ can be well ordered and no
$\lambda \in (\mu,2^\mu]$ is measurable).  
\end{theorem}

Note that before this Apter and Magidor \cite{ApMg95} had
proved the consistency of ``${\cH}(\mu)$ well
ordered, $\mu = \beth_\omega,(\forall \kappa< \mu)\text{DC}_\kappa$ 
and $\mu^+$ is measurable" so \ref{0.1} says that this consistency
result cannot be fully lifted to uncountable cofinalities.  
Generally without full choice,
a successor cardinal being not measurable is worthwhile information.

\noindent
A second theorem (\cite[\S5]{Sh:497}) was
\begin{theorem}
\label{0.2}
Assume
\mn
\begin{enumerate}
\item[$(a)$]    {\rm DC + AC}$_\kappa + \kappa$ regular uncountable
\sn
\item[$(b)$]   $\langle \mu_i:i < \kappa \rangle$ is increasing
continuous with limit $\mu,\mu > \kappa,\cH(\mu)$ is well ordered, $\mu$ 
strong limit, (we need just a somewhat weaker version, the so-called
$i < \kappa \Rightarrow Tw_{{\cD}_\kappa}(\mu_i) < \mu$).
\end{enumerate}
\mn
\Then  \ , we cannot have two regular cardinals $\theta$ such that
for some stationary $S \subseteq \kappa$, the sequence 
$\langle{\text{\rm cf\/}}(\mu^+_i):i \in S \rangle$ is constantly $\theta$.  
\end{theorem}

A dream was to prove that there is a class of regular cardinals 
from a restricted version of choice (see more \cite{Sh:497} and a 
little more in \cite{Sh:E37}).

Our original aim here is to improve those theorems.  As for \ref{0.1} we
replace ``${\cH}(\mu)$ well ordered" by ``$[\mu]^{\aleph_0}$ is
well ordered" and then by weaker statements.

We know (assuming full choice) that if, e.g., 
$\neg \exists 0^\#$ or there is 
no inner model with a measurable cardinal
then though $\langle 2^\kappa:\kappa$ regular$\rangle$ is quite
arbitrary, the size of $[\lambda]^\kappa,\lambda >> \kappa$ is strictly
controlled (by Easton forcing \cite{Est70}, and Jensen and Dodd
\cite{DJ2} respectively).  
It seemed that the situation here is parallel in some
sense; under the
restricted choice we assume, we cannot say much on the cardinality of 
${\cP}(\kappa)$ but can say something on the cardinality of 
$[\lambda]^\kappa$ for $\lambda >> \kappa$.

In the proofs we fulfill a promise from \cite[\S5]{Sh:589} about using
$J[f,D]$ from Definition \ref{0.C} instead of the nice filters used
in \cite{Sh:497} and, to some extent, in early versions of this work
which require going through inner models to prove their existence.  
This work is continued in Larson-Shelah \cite{LrSh:925} and will be
continued in \cite{Sh:1005}.  On a different
line with weak choice (say DC$_{\aleph_0} +
\text{ AC}_\mu,\mu$ fixed): see \cite{Sh:938}, \cite{Sh:955} and
\cite{Sh:F1039}.  The present work fits the thesies of \cite{Sh:g}
which in particular says:
it is better to look at $\langle
\lambda^{\aleph_0}:\lambda$ a cardinality$\rangle$ then at $\langle
2^\lambda:\lambda$ a cardinal$\rangle$.  
Here instead well ordering ${\cP}(\lambda)$ we well 
order $[\lambda]^{\aleph_0}$, this is enough for much.

A simply stated conclusion is (see \ref{c20})
\begin{conclusion}
\label{0n.2}
[DC] Assume $[\lambda]^{\aleph_0}$ is well ordered for every $\lambda$.

\noindent
1) If $2^{2^\kappa}$ is well ordered \then \, for every
$\lambda,[\lambda]^\kappa$ is well ordered.

\noindent
2) For any set $Y$, there is a derived set $Y_*$ so called
{\rm Fil}$^4_{\aleph_1}(Y)$ of power near ${\cP}({\cP}(Y))$ such that
$\Vdash_{\text{\rm Levy}(\aleph_0,Y)}$ ``for every 
$\lambda,{}^Y \lambda$ is well ordered". 
\end{conclusion}

\begin{thesis}
\label{0n.9}
1) If $\bold V \models ``\text{ZF} + \text{ DC}"$ and
``every $[\lambda]^{\aleph_0}$ is well orderable" \then \, $\bold V$
looks like the result of starting with a model of ZFC and using
$\aleph_1$-complete forcing notions like Easton forcing, Levy collapses,
and more generally, iterating of $\kappa$-complete forcing for $\kappa
> \aleph_0$.

\noindent
2) This approach is dual to investigating $\bold L[\bbR]$ - here we assume
$\omega$-sequences are understood (or weaker versions) and we try to
understand $\bold V$ (over this), there over the 
reals everything is understood.

Also though our original motivation was to look at consequences of
Ax$_4$, this was shadowed here by the try to use weaker relatives; 
see more in \cite{Sh:1005}.
\end{thesis}
\bigskip

\begin{explanation}
How do we analyze $[\mu]^\kappa$ or equivalently ${}^\kappa \mu$ here?
We use $\aleph_1$-complete filters on $\kappa$ and a well ordering of
$[\alpha]^{\aleph_0}$ for appropriate $\alpha$ or less.  We will consider
$f:\kappa \rightarrow \mu$; now for every $\aleph_1$-complete filter $D$
on $\kappa$, the ordinal $\rk_D(f)$ gives us some information on
$\alpha$, but if $A,\kappa \backslash A \in D^+$ and $f \rest A =
0_A$, then $\alpha = 0$ but we have no information on $f \rest (\kappa
\backslash A)$, then $\alpha =0$ but we have no information on $f
\rest (\kappa \backslash A)$.  Trying to correct this we consider the
ideal $J[f,D] = \{A \subseteq \kappa:A = \emptyset \mod D$ or $A \in
D^+$ but $\rk_{D +(A)}(f) > \alpha\}$, this is an $\aleph_1$-complete
ideal and so we may consider the pair $\bar D = (D_1,D_2) =
(D,\dual(J[f,D]))$.  Now $\alpha$ and the pair $\bar D$ gives more
information on $f$; they determine $f$ modulo $D_2$.  This is not enough
so we use an algebra $\cB$ on $\mu$ with no infinite decreasing
sequence of sub-algebras built using the assumption ``$[\mu]^{\aleph_0}$
is well ordered".  So there is $Z \in D_2$ such that $A = c
\ell_{\cB}(\Rang(f \rest Z))$ is $\subseteq$-minimal.

Now the triple $(D_1,D_2,Z)$ and the ordinal $\alpha$ almost
determines $f$, we need one more piece of information with domain
$\kappa:h(i) = \otp(\alpha \cap Z)$, hence an ordinal $<
\hrtg(\Rang(f))$.  So we need a bound on it which depends on the
choice of $\cB$, usually it is $\hrtg([\kappa]^{\aleph_0})$, natural
by the construction of $\cB$.

So $f \rest Z$ is uniquely determined by the ordinal $\rk_D(f)$ and
the quadruple $(D_1,D_2,Z,h)$, which belongs to a set defined from
$\kappa$, independently of $\mu$.  

Lastly, considering all such filters $D$ (recalling we are assuming
$\DC$) we can find countably many
quadruple $(D^n_1,D^n_2,Z^n,h^n)$ which together are enough as
$\bigcup\limits_{n} Z^n = \kappa$.
\end{explanation}
\bigskip

\noindent
We thank for attention and comments the audience in the advanced
seminar in Rutgers 10/2004 (particularly Arthur Apter) and advanced course
in logic in the Hebrew University 4,5/2005 and to Paul Larson and
Shimoni Garti for many corrections.

\bigskip

\subsection{Preliminaries} \label{0B}\
\bigskip

\begin{convention}
\label{0n.2.7}
We assume just $\bold V \models \text{ ZF}$ if not said otherwise.
\end{convention}

\begin{notation}
Let

\noindent
1) $\alpha,\beta,\gamma,\delta,\varp,\zeta,\xi,i,j$ denote ordinals.

\noindent
2) $\kappa,\lambda,\mu,\chi$ denotes cardinals, infinite if not said otherwise.

\noindent
3) $n,m,k,\ell$ denotes natural numbers.

\noindent
4) $D$ denotes a filter (on some set), $I,J$ denote ideals on some set. 
\end{notation}

\begin{definition}
\label{0.3}
1) \hrtg$(A) = \text{ Min}\{\alpha$: there is no 
function from $A$ onto $\alpha\}$.

\noindent
2) $\wlor(A) = \text{ Min}\{\alpha$: there is no one-to-one function from
$\alpha$ into $A$ or $\alpha= 0 \wedge A = \emptyset\}$ so
$\wlor(A) \le$ \hrtg$(A)$.
\end{definition}

\begin{definition}
\label{0.A}
1) For $D$ an $\aleph_1$-complete
filter on $Y$ and $f \in {}^Y\text{Ord}$ and $\alpha \in \text{ Ord }
\cup \{\infty\}$ we define when rk$_D(f) = \alpha$, by induction on $\alpha$:
\mn
\begin{enumerate}
\item[$\circledast$]   For $\alpha < \infty$, rk$_D(f) = \alpha$ iff
$\beta < \alpha \Rightarrow \text{ rk}_D(f) \ne \beta$ and for every $g \in
{}^Y\text{Ord}$ satisfying $g <_D f$ there is $\beta < \alpha$ such
that rk$_D(g) = \beta$.
\end{enumerate}
\mn
2) We can replace $D$ by the dual ideal.  If $f \in {}^Z\text{Ord}$
   and $Z \in D$ then we let $\rk_D(f) = \rk_{D+Z}(f \cup 0_{Y\backslash Z})$.
\end{definition}

Galvin-Hajnal \cite{GH75} use the rank for the
club filter on $\omega_1$.  This
was continued in \cite{Sh:71} where varying $D$ was extensively used.
\begin{claim}
\label{0.B}
{\rm [DC]} In Definition \ref{0.A}, {\rm rk}$_D(f)$ is always
an ordinal and if $\alpha \le \text{\rm rk}_D(f)$ \then \,
for some $g \in \prod\limits_{y \in Y} (f(y)+1)$ we 
have $\alpha = \text{\rm rk}_D(g)$, (if $\alpha <
\text{\rm rk}_D(f)$ we can add $g <_D f$; if {\rm rk}$_D(f) < \infty$
then {\rm DC} is not necessary; if {\rm rk}$_D(f) = \alpha$ this is
trivial, as we can choose $g=f$).
\end{claim}

\begin{claim}
\label{0.B1}
1) {\rm [DC]} If $D$ is an $\aleph_1$-complete
filter on $Y$ and $f \in {}^Y${\rm Ord} and $Y = 
\cup\{Y_n:n < \omega\}$ \then \, {\rm rk}$_D(f) = 
{ \text{\rm Min\/}}\{\text{\rm rk}_{D+Y_n}(f):n < \omega$ and
$Y_n \in D^+\}$, (\cite{Sh:71}).

\noindent
2) [{\rm DC + AC}$_{\alpha^*}$] If $D$ is a $\kappa$-complete filter on
$Y,\kappa$ a cardinal $> \aleph_0$ and 
$f \in {}^Y${\rm Ord} and $Y = \cup\{Y_\alpha:\alpha <
\alpha^*\},\alpha^* < \kappa$ \then \, {\rm rk}$_D(f) = 
{ \text{\rm Min\/}}\{\text{\rm rk}_{D+Y_\alpha}(f):\alpha < \alpha^*$ and
$Y_\alpha \in D^+\}$.
\end{claim}

\begin{PROOF}{\ref{0.B1}}
1) By \cite{Sh:71}, in fact, $\AC_{\aleph_0}$ suffice.

\noindent
2) By \cite{Sh:71}, in fact, $\DC$ is not necessary.
\end{PROOF}

\begin{definition}
\label{0.C}
For $Y,D,f$ as in \ref{0.A} let $J[f,D] =:
\{Z \subseteq Y:Y \backslash Z \in D$ or $Y \backslash Z \in D^+$ and 
$\rk(f)_{D+(Y \backslash Z)} > \rk_D(f)\}$. 
\end{definition}

\begin{claim}
\label{0.D}
[\rm{DC+AC}$_{< \kappa}$]  Assume $D$ 
is a $\kappa$-complete filter on $Y,\kappa > \aleph_0$.

\noindent
1) If $f \in {}^Y{\text{\rm Ord\/}}$ \then \, $J[f,D]$ is a 
$\kappa$-complete ideal on $Y$. 

\noindent
2) If $f_1,f_2 \in {}^Y\text{Ord}$ and $J = J[f_1,D] = J[f_2,D]$
then {\rm rk}$_D(f_1) < \text{\rm rk}_D(f_2) 
\Rightarrow f_1 < f_2$ mod $J$ and {\rm rk}$_D(f_1) = 
\text{\rm rk}_D(f_2) \Rightarrow f_1 = f_2$ mod $J$.
\end{claim}

\begin{PROOF}{\ref{0.D}}
Straightforward or see \cite[\S5]{Sh:589} and 
the reference there to \cite{Sh:497} (and \cite{Sh:71}).
\end{PROOF}

\begin{definition}
\label{0n.D}
1) Here $Y \le_{\text{qu}} Z$ or $|Y|
\le_{\text{qu}} |Z|$ or $|Y| \le_{\text{qu}} Z$ or $Y \le_{\text{qu}} |Z|$
means that $Y = \emptyset$ or there is a function from $Z$
(equivalently from a subset of $Z$) onto $Y$.

\noindent
2) reg$(\alpha) = \text{ Min}\{\partial:\partial \ge \alpha$ is a
regular cardinal$\}$.
\end{definition}

\begin{definition}
\label{0.E}
For a set $Y$, cardinal $\kappa$ and ordinal $\gamma$ we 
define ${\cH}_{< \kappa,\gamma}(Y)$ by induction on $\gamma$: if 
$\gamma = 0,{\cH}_{< \kappa,\gamma}(Y)=Y$, if $\gamma = \beta +1$ 
then ${\cH}_{< \kappa,\gamma}(Y) = {\cH}_{< \kappa,\beta}(Y) 
\cup \{u:u \subseteq {\cH}_{< \kappa,\beta}(Y)$ and $|u| < 
\kappa\}$ and if $\gamma$ is a limit
ordinal then ${\cH}_{< \kappa,\gamma}(Y) = \cup\{{\cH}_{<
\kappa,\beta}(Y):\beta < \gamma\}$.
\end{definition}

\begin{observation}
\label{x.5}
1) If $\lambda$ is the disjoint union of
$\langle W_z:z \in Z \rangle$ and $z \in Z \Rightarrow |W_z| <
\lambda$ and $\wlor(Z) \le \lambda$ \then \, $\lambda =
\sup\{\text{\rm otp}(W_z):z \in Z\}$ hence {\rm cf}$(\lambda) <$ \hrtg$(Z)$.

\noindent
2) If $\lambda = \cup\{W_z:z \in Z\}$ and $\wlor({\cP}(Z)) 
\le \lambda$ \then \, $\sup\{\text{\rm otp}(W_z):z \in Z\} = \lambda$.

\noindent
3) If $\lambda = \cup\{W_z:z \in Z\}$ and $|Z| < \lambda$ \then \,
$\lambda = \text{\rm sup}\{\text{\rm otp}(W_z):z \in Z\}$.

\noindent
4) If $Z \subseteq \Ord,\bar W = \langle W_\alpha:\alpha \in Z
\rangle,W_\alpha \subseteq \Ord$ and $\lambda \ge
\aleph_0,|Z|,|W_\alpha|$ for $\alpha \in Z$ \then \,
$\cup\{W_\alpha:\alpha \in Z\}$ has cardinality $\le \lambda$.
\end{observation}

\begin{PROOF}{\ref{x.5}}
1) Let $Z_1 = \{z \in Z:W_z \ne \emptyset\}$, so the
mapping $z \mapsto \text{ Min}(W_z)$ exemplifies that $Z_1$ is well 
ordered hence by the definition of $\wlor(Z_1)$ the power $|Z_1|$ 
is an aleph $< \wlor(Z_1) \le \wlor(Z)$ and by 
assumption $\wlor(Z) \le \lambda$.  Now if the desirable
conclusion fails then $\gamma^* = 
\sup(\{\text{otp}(W_z):z \in Z_1\} \cup \{|Z_1|\})$ is an 
ordinal $< \lambda$, so we can find a sequence $\langle u_\gamma:
\gamma < \gamma^* \rangle$ such that otp$(u_\gamma) \le \gamma^*,
u_\gamma \subseteq \lambda$ and $\lambda
= \cup\{u_\gamma:\gamma < \gamma^*\}$, so $\gamma^* < \lambda \le 
|\gamma^* \times \gamma^*|$, easy contradiction.

\noindent
2) For $x \subseteq Z$ let $W^*_x = \{\alpha < \lambda:(\forall z \in
Z)(\alpha \in W_z \equiv z \in x)\}$ hence $\lambda$ is the
disjoint union of $\{W^*_x:x \in {\cP}(Z) \backslash \{\emptyset\}\}$.  
So the result follows by part (1).

\noindent
3) So let $<_*$ be a well ordering of $Z$ and let $W'_z = \{\alpha \in
W_z$: if $y <_* z$ then $\alpha \notin W_y\}$, so $\langle W'_z:z
\in Z\rangle$ is a well defined sequence of pairwise disjoint sets
with union equal to $\cup\{W_z:z \in Z\} = \lambda$ and otp$(W'_z)
\le \text{ otp}(W_z)$.  Hence if $|W_z| = \lambda$ for some $z \in Z$  
the desirable conclusion is obvious, otherwise 
the result follows by part (1).

\noindent
4) Should be clear.
\end{PROOF}

\begin{definition}
\label{0.F}
1) We say that $c \ell$ is a very weak closure operation on 
$\lambda$ of character $(\mu,\kappa)$ \when \,:
\mn
\begin{enumerate}
\item[$(a)$]   $c \ell$ is a function from ${\cP}(\lambda)$ to
${\cP}(\lambda)$
\sn
\item[$(b)$]   $u \in [\lambda]^{\le \kappa} \Rightarrow |c \ell(u)| \le \mu$
\sn
\item[$(c)$]   $u \subseteq \lambda \Rightarrow u \cup\{0\} 
\subseteq c \ell(u)$, the $0$ for technical reasons.
\end{enumerate}
\mn
1A)  We say that $c \ell$ is a weak closure\footnote{so by actually
 only $c \ell \rest [\lambda]^{\le \kappa}$ count}
 operation on $\lambda$ of character $(\mu,\kappa)$ \when \, 
(a),(b),(c) above and:
\mn
\begin{enumerate}
\item[$(d)$]   $u \subseteq v \subseteq \lambda \Rightarrow u
\subseteq c \ell(u) \subseteq c \ell(v)$
\sn
\item[$(e)$]    $c \ell(u) = \cup\{c \ell(v):v \subseteq u,|v| \le \kappa\}$.
\end{enumerate}
\mn
1B) Let ``... character $(< \mu,\kappa)$ or $(\mu,< \kappa)$, or $(<
\mu,< \kappa)$" have the obvious meaning but if $\mu$ is an ordinal
not a cardinal, then ``$< \mu$" means of order type $< \mu$; similarly for
``$< \kappa$".  Let ``... character $(\mu,Y)"$ 
means ``character $(< \mu^+,< \hrtg(Y))$"

\noindent
1C) We omit the weak \when \, in addition:
\mn
\begin{enumerate}
\item[$(f)$]   $c \ell(u) = c \ell(c \ell(u))$ for $u \subseteq \lambda$.
\end{enumerate}
\mn
2) We say $\lambda$ is $f$-inaccessible when $\delta \in \lambda \cap \text{
Dom}(f) \Rightarrow f(\delta) < \lambda$.

\noindent
3) We say $c \ell:{\cP}(\lambda) \rightarrow {\cP}(\lambda)$ is well
founded \when \, for no sequence $\langle {\cU}_n:n < \omega \rangle$
of subsets of $\lambda$ do we have $c \ell({\cU}_{n+1}) \subset
{\cU}_n$ for $n < \omega$.  

\noindent
4) For $c \ell$ a partial function from ${\cP}(\alpha)$ to 
${\cP}(\alpha)$ (for simplicity assume $\alpha = \cup\{u:u \in \text{
Dom}(c \ell)\})$ let $c \ell^1_{\varepsilon,< \kappa}$ be the function 
from ${\cP}(\alpha)$ to ${\cP}(\alpha)$ defined by induction on the ordinal
$\varepsilon$ as follows:
\mn
\begin{enumerate}
\item[$(a)$]   $c \ell^1_{0,< \kappa}(u) = u$
\sn
\item[$(b)$]  $c \ell^1_{\varepsilon +1,<\kappa}(u) = \{0\} \cup c
\ell^1_{\varepsilon,< \kappa}(u) \cup \bigcup\{c \ell(v):v \subseteq c 
\ell^1_{\varepsilon,< \kappa}(u)$ and $v \in \text{ Dom}(c \ell),|v| <
\kappa\}$
\sn
\item[$(c)$]  for limit $\varepsilon$ let 
$c \ell^1_{\varepsilon,< \kappa}(u) = \cup\{c \ell^1_{\zeta,<
\kappa}(u):\zeta < \varepsilon\}$.
\end{enumerate}
\mn
4A) Instead ``$< \kappa$" we may use ``$\le \kappa$".

\noindent
5) For any function $F:[\lambda]^{\aleph_0} \rightarrow \lambda$ and
countable $u \subseteq \lambda$
we define $c \ell^2_\varepsilon(u,F)$ by induction on 
$\varepsilon \le \omega_1$
\mn
\begin{enumerate}
\item[$(a)$]  $c \ell^2_0(u,F) = u \cup \{0\}$
\sn
\item[$(b)$]  $c \ell^2_{\varepsilon +1}(u,F) = 
c \ell^2_\varp(u,F) \cup \{F(c \ell^2_\varepsilon(u,F))\}$
\sn
\item[$(c)$]  $c \ell^2_\varepsilon(u,F) = \cup\{c \ell^2_\zeta(u,F):\zeta <
\varepsilon\} \text{ when } \varepsilon \le \omega_1 \text{ is a limit
ordinal}$.
\end{enumerate}
\mn
6) For countable $u$ and $F$ as in part (5) let $c \ell^3_F(u) = c
\ell^3(u,F) := c \ell^2_{\omega_1}(u,F)$ and for any $u \subseteq \lambda$ let
$c\ell^4_F(u) := u \cup \bigcup\{c \ell^3_F(v):v \in \Dom(F)\}$.

\noindent
7) For a cardinal $\partial$ we say that  
$c \ell:{\cP}(\lambda) \rightarrow {\cP}(\lambda)$ is 
$\partial$-well founded \when \, for no $\subseteq$-decreasing
sequence $\langle {\cU}_\varepsilon:\varepsilon <
\partial\rangle$ of subsets of $\lambda$ do we have $\varepsilon <
\zeta < \partial \Rightarrow c \ell({\cU}_\zeta) 
\nsupseteq {\cU}_\varepsilon$.

\noindent
8) If $F:[\lambda]^{\le \kappa} \rightarrow \lambda$ and $u \subseteq
\lambda$ \then \, we let $c \ell_F(u) = c \ell^1_F(u)$ be the minimal
subset $v$ of $\lambda$ such that $w \in [v]^{\le \kappa}
\Rightarrow F(w) \in  v$ and $u \subseteq v$ (exists).
\end{definition}

\begin{observation}
\label{0n.F}
For $F:[\lambda]^{\aleph_0}
\rightarrow \lambda$, the operation $u \mapsto c \ell^3_F(u)$ is a very
weak closure operation of character $(\aleph_1,\aleph_0)$.
\end{observation}

\begin{remark}
\label{0.G}
So for any very weak closure operation, $\aleph_0$-well founded is a stronger
property than well founded, but if $u \subseteq \lambda \Rightarrow c
\ell(c \ell(u)) = c \ell(u)$ which is reasonable, they are equivalent.
\end{remark}

\begin{observation}
\label{0.H}
$[\alpha]^\partial$ is well ordered iff
${}^\partial \alpha$ is well ordered when $\alpha \ge \partial$.
\end{observation}

\begin{PROOF}{\ref{0.H}}
Use a pairing function on $\alpha$ for showing
$|{}^\partial \alpha| \le [\alpha]^\partial$, so $\Rightarrow$ holds.  If
${}^\partial \alpha$ is well ordered by $<_*$ map $u \in
[\alpha]^\partial$ to the $<_*$-first $f \in {}^\partial \alpha$
satisfying Rang$(f) = u$.
\end{PROOF}
\newpage

\section {Representing ${}^\kappa \lambda$} \label{1}

Here we give a simple case to illustrate what we do (see later on
improvements in the hypothesis and the conclusion).  Specifically, if $Y$
is uncountable and $[\lambda]^{\aleph_0}$ is
well ordered, \then \, the set ${}^Y \lambda$ can be analyzed modulo
countable union over few (i.e., their number depends on $Y$ but not on
$\lambda$) well ordered sets.

\begin{definition}
\label{r.1}
1) 
\mn
\begin{enumerate}
\item[$(a)$]    Fil$_{\aleph_1}(Y) = \text{ Fil}^1_{\aleph_1}(Y) =
\{D:D$ is an $\aleph_1$-complete filter on $Y\}$, so $Y$ is defined from
$D$ as $\cup\{X:X \in D\}$
\sn
\item[$(b)$]   Fil$^2_{\aleph_1}(Y) = \{(D_1,D_2):D_1 \subseteq D_2$ are
$\aleph_1$-complete filters  on $Y$, ($\emptyset \notin D_2$, of course)$\}$;
 in this context $Z \in \bar D$ means $Z \in D_2$
\sn
\item[$(c)$]   Fil$^3_{\aleph_1}(Y,\mu) = \{(D_1,D_2,h):(D_1,D_2) \in
\text{ Fil}^2_{\aleph_1}(Y)$ and $h:Y \rightarrow \alpha$ for some
$\alpha < \mu\}$, if we omit $\mu$ we mean $\mu =$ \hrtg$(Y) \cup \omega$
\sn
\item[$(d)$]    Fil$^4_{\aleph_1}(Y,\mu) = \{(D_1,D_2,h,Z):
(D_1,D_2,h) \in \text{ Fil}^3_{\aleph_1}(Y,\mu),Z \in D_2\}$; omitting
$\mu$ means as above. 
\end{enumerate}
\mn
2) For ${\gy} \in \text{ Fil}^4_{\aleph_1}(Y,\mu)$ let 
$Y = Y^{[\gy]} = Y[\gy]$ and $\gy = (D^{\gy}_1,D^{\gy}_2,h^{\gy},Z^{\gy}) 
= (D_1[{\gy}],D_2[{\gy}],h[{\gy}],Z[{\gy}])$; similarly for the
others and let $D^{\gy} = D[\gy]$ be $D^{\gy}_1 + Z^{\gy}$.

\noindent
3) We can replace $\aleph_1$ by any $\kappa > \aleph_1$ (the results
can be generalized easily assuming DC + AC$_{< \kappa}$, used in \S2).
\end{definition}

\begin{theorem}
\label{r.2}
[$\DC$] Assume $[\lambda]^{\aleph_0}$ is well ordered.

\Then \, we can find a sequence $\langle {\cF}_{\gy}:{\gy} 
\in { \text{\rm Fil\/}}^4_{\aleph_1}(Y)\rangle$ satisfying
\mn
\begin{enumerate}
\item[$(\alpha)$]  ${\cF}_{\gy} \subseteq {}^{Z[\gy]}\lambda$
\sn
\item[$(\beta)$]    ${\cF}_{\gy}$ is a well ordered set by
$f_1 <_{\frak y} f_2 \Leftrightarrow \text{\rm rk}_{D[{\gy}]}
(f_1) < \text{\rm rk}_{D[{\gy}]} (f_2)$ so $f \mapsto 
\text{\rm rk}_{D[{\gy}]}(f)$ is a one-to-one mapping from
${\cF}_{\gy}$ into the ordinals
\sn
\item[$(\gamma)$]   if $f \in {}^Y \lambda$ \then \, we can find a sequence
$\langle {\gy}_n:n < \omega \rangle$ with ${\gy}_n
\in { \text{\rm Fil\/}}^4_{\aleph_1}(Y)$ such that 
$n < \omega \Rightarrow f \restriction Z^{{\gy}_n} 
\in {\cF}_{{\gy}_n}$ and $\cup\{Z^{{\gy}_n}:n < \omega\} = Y$.
\end{enumerate}
\end{theorem}

\noindent
An immediate consequence of \ref{r.2} is
\begin{conclusion}
\label{r.1x}
1) [DC + ${}^\omega \alpha$ is well-orderable for every ordinal $\alpha$].

For any set $Y$ and cardinal $\lambda$ there is a sequence
$\langle{\cF}_{\bar{\gx}}:\bar{\gx} \in
{}^\omega(\text{\rm Fil}^4_{\aleph_1}(Y))\rangle$ such that
\mn
\begin{enumerate}
\item[$(a)$]   ${}^Y \lambda = \cup\{{\cF}_{\bar{\gx}}:
\bar{\gx} \in {}^\omega(\text{\rm Fil}^4_{\aleph_1}(Y))\}$
\sn
\item[$(b)$]   ${\cF}_{\bar{\gx}}$ is well orderable for each
$\bar{\frak x} \in {}^\omega(\text{\rm Fil}^4_{\aleph_1}(Y))$
\sn
\item[$(b)^+$]   moreover, uniformly, i.e., there is a sequence
$\langle <_{\bar{\gx}}:\bar{\gx} \in
{}^\omega(\text{\rm Fil}^4_{\aleph_1}(Y)\rangle$ such that $<_{\bar{\gx}}$
is a well order of ${\cF}_{\bar{\gx}}$
\sn
\item[$(c)$]   there is a function $F$ with domain ${\cP}({}^Y
  \lambda) \backslash \{\emptyset\}$ such that: if $S \subseteq {}^Y
  \lambda$ is non-empty 
then $F(S)$ is a non-empty subset of $S$ of power 
$\le_{\qu} {}^\omega(\text{\rm Fil}^4_{\aleph_1}(Y)))$ recalling Definition
\ref{0n.D}.  In fact, some ordinal $\alpha(*)$ and $\bar u$ we have:
\sn
\begin{enumerate}
\item[$(\alpha)$]  $\bar u = \langle {\cU}_\alpha:
\alpha < \alpha(*)\rangle$ is a partition of ${}^Y \lambda$
\sn
\item[$(\beta)$]  if $S \subseteq {}^Y \lambda$ then 
$F(S) = {\cU}_{f(S)} \cap S$ where 
$f(S) = \text{\rm Min}\{\alpha:{\cU}_\alpha \cap S
\ne \emptyset\}$
\sn
\item[$(\gamma)$]  if $\alpha < \alpha(*)$ then $|{\cU}_\alpha| <$ 
\hrtg$({}^\omega(\text{\rm Fil}^4_{\aleph_1}(Y)))$. 
\end{enumerate}
\end{enumerate}
\mn
2) [$\DC$] For any $Y,\lambda$ above, if 
$[\alpha(*)]^{\aleph_0}$ is well ordered where $\alpha(*) =
\cup\{\text{\rm rk}_D(f) + 1:f \in {}^Y \lambda$ and $D \in 
\text{\rm Fil}^1_{\aleph_1}(Y)\}$ \then \, ${}^Y \lambda$ satisfies the
conclusion of part (1).
\end{conclusion}

\begin{remark}
So clause (c) of \ref{r.1x}(1) is a weak form of choice.
\end{remark}

\begin{PROOF}{\ref{r.1x}}
\underline{Proof of \ref{r.1x}}  
1) Let $\langle {\cF}_{\gy}:{\gy} \in 
\text{ Fil}^4_{\aleph_1}(Y)\rangle$ be as in \ref{r.2}.

For each $\bar{\gx} \in {}^\omega(\text{Fil}^4_{\aleph_1}(Y))$ (so
$\bar{\gx} = \langle {\frak x}_n:n < \omega \rangle$) let

\begin{equation*}
\begin{array}{clcr}
{\cF}'_{\bar{\gx}} = \{f: &f \text{ is a function from } Y
\text{ to } \lambda \text{ such that} \\
  &n < \omega \Rightarrow f \restriction Z^{{\gx}_n} \in 
{\cF}_{{\gx}_n} \text{ and }  Y = \cup\{Z^{{\gx}_n}:n < \omega\}\}.
\end{array}
\end{equation*}

\mn
Now
\mn
\begin{enumerate}
\item[$(*)_1$]   ${}^Y \lambda = \cup\{{\cF}'_{\bar{\gx}}:
\bar{\gx} \in {}^\omega(\text{Fil}^4_{\aleph_1}(Y))\}$.
\end{enumerate}
\mn
[Why?  By clause $(\gamma)$ of \ref{r.2}.]

Let $\alpha(*) = \cup\{\text{rk}_D(f) +1:f \in {}^Y \lambda$ and $D \in
\text{ Fil}^1_{\aleph_1}(Y)\}$.  For $\bar{\gx} \in
{}^\omega(\text{Fil}^4_{\aleph_1}(Y))$ we define the function
$G_{\bar{\gx}}:{\cF}'_{\bar{\gx}} \rightarrow {}^\omega
\alpha(*)$ by $G_{\bar{\gx}}(f) = \langle 
\text{\rm rk}_{D_1[{\gx}_n]}(f):n  < \omega \rangle$.

Next
\mn
\begin{enumerate}
\item[$(*)_2$]   $(\alpha) \quad \bar G = \langle 
G_{\bar{\gx}}:\bar{\gx} \in {}^\omega(\text{Fil}^4_{\aleph_1}(Y))\rangle$
exists
\sn
\item[${{}}$]   $(\beta) \quad G_{\bar{\gx}}$ is a function from
${\cF}'_{\bar{\gx}}$ to ${}^\omega \alpha(*)$
\sn
\item[${{}}$]   $(\gamma) \quad G_{\bar{\gx}}$ is one to one.
\end{enumerate}
\mn
[Should be clear, e.g. for $(*)_2(\gamma)$ read the definition of
$\cF'_{\gx'}$ and clause $(\beta)$ of Theorem \ref{r.2}.]

Let $<_*$ be a well ordering of ${}^\omega \alpha(*)$ and for
$\bar{\gx} \in {}^\omega(\text{Fil}^4_{\aleph_1}(Y))$ let
$<_{\bar{\gx}}$ be the following two place relation on ${\cF}'_{\bar{\gx}}$:
\mn
\begin{enumerate}
\item[$(*)_3$]   $f_1 <_{\bar{\gx}} f_2$ \Iff \, 
$G_{\bar{\gx}}(f_1) <_* G_{\bar{\gx}}(f_2)$.
\end{enumerate}
\mn
Obviously
\mn
\begin{enumerate}
\item[$(*)_4$]   $(\alpha) \quad \langle <_{\bar{\gx}}:\bar{\gx} 
\in {}^\omega(\text{Fil}^4_{\aleph_1}(Y))\rangle$ exists
\sn
\item[${{}}$]   $(\beta) \quad <_{\bar{\gx}}$ is a well ordering
of ${\cF}'_{\bar{\gx}}$.
\end{enumerate}
\mn
By $(*)_1 + (*)_4$ we have proved clauses (a),(b),(b)$^+$ of the
conclusion.  Now clause (c) follows: for non-empty 
$S \subseteq {}^Y \lambda$, let $f(S)$ be
min$\{\text{otp}(\{g:g <_{\bar{\gy}} f\},<_{\bar{\gy}}):
\bar{\gy} \in {}^\omega(\text{Fil}^4_{\aleph_1}(Y))$ and $f
\in {\cF}'_{\bar{\frak y}} \cap S\}$.  Also for any ordinal $\gamma$ 
let ${\cU}^1_\gamma := \{f$: for some $\bar{\gy} 
\in {}^\omega(\text{Fil}^4_{\aleph_1}(Y))$ we have
$\gamma = \text{ otp}(\{g:g <_{\bar{\gy}} f\},<_{\bar{\gy}})\}$
and ${\cU}_\gamma = {\cU}^1_\gamma \backslash \cup
\bigcup\{{\cU}^1_\beta:\beta < \gamma\}$.

Lastly, we let $F(S) = {\cU}_{f(S)} \cap S$.  Now check. 

\noindent
2) Similarly.
\end{PROOF}

\begin{PROOF}{\ref{r.2}}
\underline{Proof of Theorem \ref{r.2}}  
First
\mn
\begin{enumerate}
\item[$\circledast_1$]    there are a cardinal $\mu$ and a sequence 
$\bar u = \langle u_\alpha:\alpha < \mu \rangle$ listing 
$[\lambda]^{\aleph_0}$.
\end{enumerate}
\mn
[Why?  By the assumption.]

Second, we can deduce
\mn
\begin{enumerate}
\item[$\circledast_2$]    there are $\mu_1 \le \mu$ and a sequence
$\bar u = \langle u_\alpha:\alpha < \mu_1 \rangle$ such that:
\begin{enumerate}
\item[$(a)$]   $u_\alpha \in [\lambda]^{\aleph_0}$
\sn
\item[$(b)$]   if $u \in [\lambda]^{\le \aleph_0}$ then for some
finite $w \subseteq \mu_1,u \subseteq \cup\{u_\beta:\beta \in w\}$
\sn
\item[$(c)$]    $u_\alpha$ is not included in $u_{\alpha_0} \cup
\ldots \cup u_{\alpha_{n-1}}$ when $n <
\omega,\alpha_0,\dotsc,\alpha_{n-1} < \alpha$.
\end{enumerate}
\end{enumerate}
\mn
[Why?  Let $\bar u^0$ be of the form $\langle u_\alpha:\alpha <
\alpha^*\rangle$ such that $(a) + (b)$ holds and $\ell g(\bar
u^0)$ is minimal; it is well defined and $\ell g(\bar u^0) \le \mu$ by
$\circledast_1$.  Let $W = \{\alpha < \ell g(\bar u^0):u^0_\alpha \nsubseteq 
\cup\{u^0_\beta:\beta \in w\}$ when $w \subseteq \alpha$ 
is finite$\}$.  Let $\mu_1 = |W|$
and let $f:\mu_1 \rightarrow W$ be one-to-one onto, let 
$u_\alpha = u^0_{f(\alpha)}$ so $\langle u_\alpha:\alpha <
\mu_1 \rangle$ satisfies $(a) + (b)$ and $\mu_1 = |W| \le 
\ell g(\bar u^0)$.  So by the choice of $\bar u^0$ we have 
$\ell g(\bar u^0) = \mu_1$.  So we can choose $f$ such that it is
increasing hence $\bar u$ is as required.]
\mn
\begin{enumerate}
\item[$\circledast_3$]    we can define $\bold n:[\lambda]^{\le
\aleph_0} \rightarrow \omega$ and partial functions
$F_\ell:[\lambda]^{\le \aleph_0} \rightarrow
\mu_1$ for $\ell < \omega$ (so $\langle F_\ell:\ell < \omega\rangle$ exists)
 as follows:
\begin{enumerate}
\item[$(a)$]  $u$ infinite $\Rightarrow F_0(u) = \text{ Min}\{\alpha$: for some
finite $w \subseteq \alpha,u \subseteq u_\alpha \cup
\bigcup\{u_\beta:\beta \in w\}$ mod finite$\}$
\sn
\item[$(b)$]  $u$ finite $\Rightarrow F_0(u)$ undefined
\sn
\item[$(c)$]  $F_{\ell +1}(u) := F_0(u \backslash (u_{F_0(u)} \cup \ldots \cup
u_{F_\ell(u)}))$ for $\ell < \omega$ when $F_\ell(u)$ is defined
\sn
\item[$(d)$]  $\bold n(u) := \text{ Min}\{\ell:F_\ell(u) \text{
undefined}\}$.
\end{enumerate}
\end{enumerate}
\mn
Then
\mn
\begin{enumerate}
\item[$\circledast_4$]  
\begin{enumerate}
\item[(a)]  $F_{\ell +1}(u) < F_\ell(u) < \mu_1$ when they are well defined
\sn
\item[(b)]  $\bold n(u)$ is a well defined natural number and $u
\backslash \cup\{u_{F_\ell(u)}:\ell < \bold n(u)\}$ is finite and $k <
\bold n(u) \Rightarrow (u \backslash \cup\{u_{F_\ell(u)}:\ell < k\})
\cap u_{F_k(u)}$ is infinite
\sn
\item[(c)]  if $u_1,u_2 \in [\lambda]^{\aleph_0},u_1
\subseteq u_2$ and $u_2 \backslash u_1$ is finite then $F_\ell(u_1) =
F_\ell(u_2)$ for $\ell < \bold n(u_1)$ and $\bold n(u_1) = \bold n(u_2)$
\end{enumerate}
\sn
\item[$\circledast_5$]   define $F_*:[\lambda]^{\aleph_0}
\rightarrow \lambda$ by $F_*(u) = \text{ Min}(\cup\{u_{F_\ell(u)}:\ell <
\bold n(u)\} \cup \{0\} \backslash u)$ if well defined, zero otherwise
\newline
[Note: the reader may wonder: if you add $\{0\}$ then $\Min(-) = 0$ in
all cases.  However, if $0 \in u$ then by ``$\backslash u$", zero does
not belong to the set from which we choose a minimal ordinal.]
\sn
\item[$\circledast_6$]   if $u \in [\lambda]^{\aleph_0}$ then
\sn
\begin{enumerate}
\item[$(\alpha)$]  $c \ell^3(u,F_*) = c \ell^3_{F_*}(u)$ is 
$F'(u) := u \cup \bigcup\{u_{F_\ell(u)}:\ell < \bold n(u)\} \cup
\{0\}$
\sn
\item[$(\beta)$]  $c \ell^3_{F_*}(u) = c \ell^2_{\varepsilon(u)}(F)$
for some $\varepsilon(u) < \omega_1$
\sn
\item[$(\gamma)$]   there is $\bar F = \langle
F'_\varepsilon:\varepsilon < \omega_1\rangle$ such that: for every $u
\in [\lambda]^{\aleph_0},c \ell^3_{F_*}(u) =
\{F'_\varepsilon(u):\varepsilon < \varepsilon(u)\}$ and
$F'_\varepsilon(u)=0$ if $\varepsilon \in [\varepsilon(u),\omega_1)$
\sn
\item[$(\delta)$]  in fact $F'_\varepsilon(u)$ is the $\varepsilon$-th
member of $c \ell^3_{F_*}(u)$ if $\varepsilon < \varepsilon(u)$.
\end{enumerate}
\end{enumerate}
\mn
[Why?  Define $w^\varepsilon_u$ by induction on $\varepsilon$ by
$w^0_u=u,w^{\varepsilon +1}_u = w^\varepsilon_u \cup
\{F_*(w^\varepsilon_u)\}$ and for limit ordinal $\varepsilon$ we let
$w^\varepsilon_u =
\cup\{w^\zeta_u:\zeta < \varepsilon\}$.  We can prove by induction on
$\varepsilon$ that $w^\varepsilon_u \subseteq F'(u)$ which is
countable.  The partial function $g$ with domain $F'(u) \backslash u$
to Ord, $g(\alpha) = \text{ Min}\{\varepsilon:\alpha \in
w^{\varepsilon +1}_u\}$ is one to one onto an ordinal call it
$\varepsilon(*)$, so $w^{\varepsilon(*)}_u \subseteq F'(u)$ and if
they are not equal that $F_*(w^{\varepsilon(*)}_u) \in F'(u)
\backslash w^{\varepsilon(*)}_u$ hence $w^{\varepsilon(*)}_u
\subsetneqq w^{\varepsilon(*)+1}_u$ contradicting the choice of
$\varepsilon(*)$.  So clause $(\alpha)$ holds.  In fact, $c
\ell^3(u,F_*) = w^{\varepsilon(*)}_u$ and clause $(\beta)$ holds.
CLauses $(\gamma),(\delta)$ should be clear.]
\mn
\begin{enumerate}
\item[$\circledast_7$]   there is no sequence $\langle {\cU}_n:n
< \omega \rangle$ such that:
\sn
\begin{enumerate}
\item[$(a)$]   ${\cU}_{n+1} \subseteq {\cU}_n \subset \lambda$
\sn
\item[$(b)$]   ${\cU}_n$ is closed under $F_*$, i.e. $u \in
[\cU_n]^{\aleph_0} \Rightarrow F_*(u) \in \cU_n$
\sn
\item[$(c)$]   ${\cU}_{n+1} \ne {\cU}_n$.
\end{enumerate}
\end{enumerate}
\mn
[Why?  Assume toward contradiction that $\langle {\cU}_n:n <
\omega\rangle$ satisfies clauses (a),(b),(c).
Let $\alpha_n = \text{ Min}({\cU}_n \backslash {\cU}_{n+1})$ for $n <
\omega$ hence the sequence $\bar\alpha
= \langle \alpha_n:n < \omega\rangle$ is well defined 
with no repetitions and let
$\beta_{m,\ell} := F_\ell(\{\alpha_n:n \ge m\})$ for $m < \omega$ and
$\ell < \bold n_m := \bold n(\{\alpha_n:
n \in [m,\omega)\})$.  As $\bar\alpha$ is with no repetition, $\bold
n_m  >0$ and by $\circledast_4(c)$ clearly $\bold n_m = \bold n_0$ for
$m < \omega$ and $\beta_{m,\ell} = \beta_{m,0}$ for $m < \omega,\ell <
\bold n_0$.  So letting $v_m = \cup\{u_{F_\ell(\{\alpha_n:n \in
[m,\omega)\})}:\ell < \bold n_m\}$, it does not depend on $m$ so $v_m =
v_0$, and by the choice of $F_*$, as $\{\alpha_n:n \in [m,\omega)\}
\subseteq \cU_m$ and $\cU_m$ is closed under $F_*$ clearly $v_m
\subseteq \cU_m$.  Together $v_0 = v_m \subseteq \cU_m$ so $v_0
\subseteq \cap\{\cU_m:m < \omega\}$.  Also, by the definition of the
$F_\ell$'s, $\{\alpha_n:n < \omega\} \backslash v_0$ is finite so for
some $k < \omega,\{\alpha_m:n \in [k,\omega)\} \subseteq v_0$ but $v_0
\subseteq \cU_{k+1}$ contradicting the choice of $\alpha_k$.]

Moreover, recalling Definition \ref{0.F}(6):
\mn
\begin{enumerate}
\item[$\circledast'_7$]   there is no sequence $\langle {\cU}_n:n
< \omega\rangle$ such that
\begin{enumerate}
\item[$(a)$]   ${\cU}_{n+1} \subseteq {\cU}_n \subseteq \lambda$
\sn
\item[$(b)$]   $c \ell^4_{F_*}({\cU}_{n+1}) \backslash {\cU}_n \ne \emptyset$.
\end{enumerate}
\end{enumerate}
\mn
[Why?  As above but letting $\alpha_n = \text{ Min}({\cU}_n \backslash c
\ell^3_{F_*}({\cU}_n))$.]

Now we define for $(D_1,D_2,h,Z) \in \text{ Fil}^4_{\aleph_1}(Y)$ and
ordinal $\alpha$ the following,
recalling Definition \ref{0.F}(6) for clauses (e),(f): 
\mn
\begin{enumerate}
\item[$\circledast_8$]   ${\cF}_{(D_1D_2,h,Z),\alpha} =: \{f:(a)
\quad f$ is a function from $Z$ to $\lambda$

\hskip84pt $(b) \quad \text{\rm rk}_{D_1 + Z}(f \cup 0_{(Y \backslash Z)}) 
= \alpha$

\hskip84pt $(c) \quad D_2 = \{Y \backslash X:
X \subseteq Y$ satisfies $X = \emptyset$ mod $D_1$ 

\hskip104pt or $X \in D_1^+$ and 
$\text{\rm rk}_{D_1 + X}(f \cup 0_{(Y \backslash Z)}) 
> \alpha$ 

\hskip104pt $\text{ that is rk}_{D_1+X}(f) > \alpha\}$

\hskip84pt  $(d) \quad Z \in D_2$, really follows

\hskip84pt  $(e) \quad$ if $Z' \subseteq Z \wedge Z' \in D_2$ then 

\hskip104pt $c \ell^3_{F_*}(\text{Rang}(f \restriction Z')) 
= c \ell^3_{F_*}(\text{Rang}(f))$

\hskip80pt  $(f) \quad y \in Z \Rightarrow f(y) =$ the $h(y)$-th
member of $c \ell^3_{F_*}(\text{Rang}(f))\}$.
\end{enumerate}
\mn
So we have:
\mn
\begin{enumerate}
\item[$\circledast_9$]   ${\cF}_{(D_1,D_2,h,Z),\alpha}$ has at most one
member; call it $f_{(D_1,D_2,h,Z),\alpha}$ (when defined;
pedantically we should write $f_{(D_1,D_2,h,Z),c \ell,\alpha}$)
\sn
\item[$\circledast_{10}$]   ${\cF}_{(D_1,D_2,h,Z)} =:
\cup\{{\cF}_{(D_1,D_2,h,Z),\alpha}:\alpha$ an ordinal$\}$ is a well
ordered set.
\end{enumerate}
\mn
[Why?  Define $<_{(D_1,D_2,h,Z)}$ by the $\alpha$'s, i.e. $f^1 < f^2$
\Iff \, there are $\alpha_1 < \alpha_2$ such that $f^\ell =
f_{(D_1,D_2,h,2),\alpha_\ell}$ for $\ell=1,2$.]
\mn
\begin{enumerate}
\item[$\circledast_{11}$]   if $f:Y \rightarrow \lambda$ and $Z
\subseteq Y$ then the set $\text{Rang}(f \restriction Z)$ 
has cardinality $< \hrtg(Z)$.
\end{enumerate}
\mn
[Why?  By the definition of $\hrtg(-)$ this should be clear.]
\mn
\begin{enumerate}
\item[$\circledast_{12}$]   if $f:Z \rightarrow \lambda$ and $Z
\subseteq Y$ \then \, $c \ell^4_{F_*}(\text{Rang}(f)) \subseteq
\lambda$ has cardinality $< \hrtg([Z]^{\aleph_0})$ or is finite.
\end{enumerate}
\mn
Why?  If Rang$(f)$ is countable more holds by \ref{0n.F}.
Otherwise, by $\circledast_6(\beta)$ recallng Definition \ref{0.F}(6) we have
$c \ell^4_{F_*}(\text{Rang}(f)) = \text{ Rang}(f)
\cup \{F'_\varepsilon(u):u \in [\text{Rang}(f)]^{\aleph_0}$ and
$\varepsilon < \omega_1\}$.

Let $\alpha(*)$ be minimal such that Rang$(f) \cap \alpha(*)$ has
order type $\omega_1$.  Let $h_1,h_2:\omega_1 \rightarrow \omega_1$ be
such that $h_\ell(\varepsilon) < \text{ max}\{\varepsilon,1\}$ and for
every $\varepsilon_1,\varepsilon_2 < \omega_1$ 
there is $\zeta \in [\varepsilon_1
+ \varepsilon_2 +1,\omega_1)$ such that $h_\ell(\zeta) =
\varepsilon_\ell$ for $\ell=1,2$.  Define $F:[Z]^{\aleph_0}
\rightarrow \lambda$ as follows: if $u \in
[\text{Rang}(f)]^{\aleph_0}$, let $\varepsilon_\ell(u) =
h_\ell(\text{otp}(u \cap \alpha(*))$ for $\ell=1,2$ and $F(u) 
= F'_{\varp_2(u)}(\{\alpha \in u$: if $\alpha < \alpha(*)$ then otp$(u \cap
\alpha) < \varepsilon_1(u)\})$.

Now
\mn
\begin{enumerate}
\item[$\bullet_1$]  if $u \in [\text{Rang}(f)]^{\aleph_0}$ then $F(u)$
is $F_\varepsilon(v)$ for some $v \in [Z]^{\aleph_0}$ and $\varepsilon
< \omega_1$.
\end{enumerate}
\mn
[Why?  As $F(u) \in \text{ Rang}(F'_{\varepsilon_2(u)} \rest
[\text{Rang}(f)]^{\aleph_0})$]
\mn
\begin{enumerate}
\item[$\bullet_2$]   $\{F(u):u \in [\text{Rang}(f)]^{\aleph_0}\}
\subseteq c \ell^4_{F_*}(\text{Rang}(f))$.
\end{enumerate}
\mn
[Why?  By $\bullet_1$ recalling $\circledast_6$.]
\mn
\begin{enumerate}
\item[$\bullet_3$]   if $u \in [\text{Rang}(f)]^{\aleph_0}$ and
$\varepsilon < \omega_1$ then $F'_\varepsilon(u)$ is $F(u)$ for some
$v \in [\text{Rang}(f)]^{\aleph_0}$.
\end{enumerate}
\mn
[Why? Let $\varepsilon_1 = \text{ otp}(u \cap \alpha(*)),\varepsilon
_2 = \varepsilon$; now let $\zeta < \omega_1$ be such that $h_\ell(\zeta) =
\varepsilon _\ell$ for $\ell=1,2$.  Let $v = u \cup \{\alpha:\alpha
\in \text{ Rang}(f) \cap \alpha(*)$ and $\alpha \ge \sup(u \cap
\alpha(*))+1$ and otp$(\text{Rang}(f) \cap \alpha \backslash (\sup(u
\cap \alpha(*)+1)) < (\zeta - \varepsilon_1))\}$.]

So $F(u) = F'_\varepsilon(u)$.  By $\bullet_2 + \bullet_3$ we can conclude:
\mn
\begin{enumerate}
\item[$\bullet_4$]   in $\bullet_2$ we have equality.
\end{enumerate}
\mn
Together $c \ell^4_{F_*}(\text{Rang}(f)) = \{F(u):u \in
[\text{Rang}(f)]^{\aleph_0}\} \cup \text{ Rang}(f)$ so it is the union
of two sets; by the definition of $\hrtg(-)$ 
the first is of cardinality $< \hrtg([Z]^{\aleph_0})$
and the second is of cardinality $< \hrtg[Z]$, so we are easily done
proving $\circledast_{12}$
\mn
\begin{enumerate}
\item[$\circledast_{13}$]  if $f:Y \rightarrow \lambda$ then for
some sequence $\langle (\gy_n,\alpha_n):n < \omega\rangle$ we have
${\gy}_n \in \text{ Fil}^4_{\aleph_1}(Y)$ and $\alpha_n \in 
\text{ Ord}$ for $n < \omega$ and $f = \cup\{f_{{\gy}_n,
\alpha_n}:n < \omega\}$.
\end{enumerate}
\mn
[Why?  Let 

\begin{equation*}
\begin{array}{clcr}
{\cI}^0_f = \{Z \subseteq Y:&\text{ for some } {\gy} 
\in \text{ Fil}^4_{\aleph_1}(Y) \text{ satisfying } Z^{\gy} = Z \\
  & \text{ and ordinal } \alpha,
f_{{\gy},\alpha} \text{ is well defined and equal to } f \restriction
Z\}
\end{array}
\end{equation*}

\[
{\cI}_f = \{Z \subseteq Y:Z \text{ is included in a countable union
of members of } {\cI}^0_f\}.
\]

\mn
So recalling we are assuming $\DC$ it is enough to show that $Y \in {\cI}_f$.

Toward contradiction assume not.  Let $D_1 = \{Y \backslash Z:Z \in
{\cI}_f\}$, clearly it belongs to Fil$_{\aleph_1}(Y)$, noting
that $\emptyset \in \cI_f$.  So $\alpha(*) :=
\text{\rm rk}_{D_1}(f)$ is well defined (by \ref{0.B}) recalling that only 
DC = DC$_{\aleph_0}$ is needed.

Let

\[
D_2 = \{X \subseteq Y:X \in D_1 \text{ or } 
\rk_{D_1+(Y \backslash X)}(f) > \alpha(*)\}.
\]

\mn
By \ref{0.C} + \ref{0.D} clearly 
$D_2$ is an $\aleph_1$-complete filter on $Y$ extending $D_1$.

Now we try to choose $Z_n \in D_2$ for $n < \omega$ 
such that $Z_{n+1} \subseteq Z_n$ and 
$c \ell^4_{F_*}(\text{Rang}(f \restriction Z_{n+1}))$ does not include
Rang$(f \restriction Z_n)$.

For $n=0,Z_0 = Y$ is O.K.

By $\circledast'_7$ we cannot have such $\omega$-sequence $\langle Z_n:n
< \omega\rangle$; so by DC for some (unique) $n=n(*),Z_n$ 
is chosen but not $Z_{n+1}$.  

Let $h:Z_n \rightarrow$ \hrtg$([Y]^{\aleph_0}) \cup \omega_1$ be:

\[
h(y) = \text{ otp}(f(y) \cap c \ell^4_{F_*}(\text{Rang}(f \restriction
Z_n))).
\]

\mn
Now $h$ is well defined by $\circledast_{12}$.  Easily

\[
f \restriction Z_n \in {\cF}_{(D_1 + Z_n,D_2,h,Z_n),\alpha(*)}
\]

\mn
hence $Z_n \in {\cI}^0_f \subseteq {\cI}_f$, 
contradiction to $Z_n \in D_2,D_1 \subseteq D_2$.

So we are done proving $\circledast_{13}$.]

Now clause $(\beta)$ of the conclusion holds by the
definition of ${\cF}_{\gy}$, clause $(\alpha)$ holds by
$\circledast_{10}$ recalling $\circledast_8,\circledast_9$ 
and clause $(\gamma)$ holds by $\circledast_{12}$.
\end{PROOF}

\begin{remark}
\label{r.2y}
We can improve \ref{r.2} in some way by weakening the demands on $\bar u$.  

We may replace the assumption ``$[\lambda]^{\aleph_0}$ is well ordered" by:
\mn
\begin{enumerate}
\item[$(*)$]   there is $\langle u_\alpha:\alpha < \alpha^* \rangle$,
a sequence of members of $[\lambda]^{\aleph_0}$ such that $(\forall
u \in [\lambda]^{\aleph_0})(\exists \alpha)(u \cap u_\alpha$
infinite).
\end{enumerate}
\mn
[Why?  We define $F_\varp:[\lambda]^{\aleph_0} \rightarrow \alpha^*$
by induction on $\varp < \omega_1$ by $F_\varp(v) := 
\text{ Min}\{\alpha < \alpha^*:(v \backslash v \cup \{F_*(v):\zeta <
\varp\}) \cap u_\alpha$ infinite$\}$ if well defined and let
$F:[\lambda]^{\aleph_0} \rightarrow [\lambda]^{\aleph_0}$ be defined
by $F(v) = \cup\{F_\varp(v):\varp < \omega_1,F_\varp(v)$ well
defined$\}$.

Lastly, let $F_*(u) = \text{ min}(F(u) \backslash u)$.] 
\end{remark}

\begin{observation}
\label{1p.4}
1) The power of {\rm Fil}$^4_{\aleph_1}(Y,\mu)$ is
smaller or equal to the power of the set 
$({\cP}({\cP}(Y)))^2 \times \cP(Y) \times
\mu^{|Y|}$; if $\aleph_0 \le |Y|$ this is equal to the power of
$\cP(\cP(Y)) \times {}^Y \mu$.

\noindent
2) The power of 
{\rm Fil}$^4_{\aleph_1}(Y)$ is smaller or equal to the power of 
the set $({\cP}({\cP}(Y)))^2 \times {\cP}(Y) \times \cup\{{}^Y \alpha:\alpha <
\hrtg([Y]^{\aleph_0})\}$.

\noindent
3) In part (2), if $\aleph_0 \le |Y|$ this is equal to $|\cP(\cP(Y))|
\times \cup\{{}^Y \alpha:\alpha < \hrtg([Y]^{\aleph_0})\}$;
 also $\alpha < \hrtg([Y]^{\aleph_0}) \Rightarrow
|\cP(\cP(Y)) \times {}^Y \alpha| = |\cP(\cP(Y))|$ and
$|\Fil^4_{\aleph_1}(Y)| \le_{\qu} \cP(\cP(Y \times Y))$. 
\end{observation}

\begin{remark}
1) As we are assuming DC, the case $\aleph_0 \nleq |Y|$ means that $Y$
is finite, so degenerated.  Also if $|Y| = \aleph_0$ then
$\Fil^1_{\aleph_1}(Y) = \{\{Z \subseteq Y:Z \supseteq X\}:X \subseteq
Y\}$ hence $|\Fil^1_{\aleph_1}(Y)| = \cP(Y)|$ hence
$\FIL^4_{\aleph_1}(Y,\mu)$ has the same power as $\cP(Y) \times
{}^\omega \mu$ again this is a dull case.
\end{remark}

\begin{PROOF}{\ref{1p.4}} 
1) Reading the definition of Fil$^4_{\aleph_1}(Y,\mu)$ clearly its
   power is $\le$ the power of $\cP(\cP(Y)) \times \cP(\cP(Y)) 
\times \cP(Y) \times \mu^{|Y|}$.  
If $\aleph_0 \le |Y|$ then $|\cP(\cP(Y)) \times
\cP(Y)| \le |\cP(\cP(Y)) \times \cP(\cP(Y))| 
= 2^{|\cP(Y))} \times 2^{|\cP(Y)|} \le 2^{|\cP(Y)|+|\cP(Y)|} =
2^{|\cP(Y)|} = |\cP(\cP(Y))| \le |\cP(\cP(Y)) \times \cP(Y) \times 
\mu^{|Y|}|$ as $\cP(Y) + \cP(Y) = 2^{|Y|} \times 2 = 
2^{|Y|+1} = 2^{|Y|}$; so the second conclusion follows.

\noindent
2) Read the definitions.

\noindent
3) If $\alpha < \hrtg([Y]^{\aleph_0})$ then let $f$ be a function
   from $[Y]^{\aleph_0}$ onto $\alpha$ and for $\beta < \alpha$ let
   $A_{f,\beta} = \{u \in [Y]^{\aleph_0}:f(u) < \beta\}$.  So $\beta
   \mapsto A_{f,\beta}$ is a one-to-one function from $\alpha$ onto
$\{A_{f,\gamma}:\gamma < \alpha\} \subseteq \cP(\cP(Y))$ so $|{}^Y
 \alpha| \le \cP(\cP(Y))$ and $\cP(\cP(Y)) \times |{}^Y \alpha| \le
\cP(\cP(Y)) \times \cP(\cP(Y)) \le 2^{|\cP(Y)|+|\cP)Y)|} = 2^{|\cP(Y)|}$.  
Better, for $f$ a function from
$[Y]^{\aleph_0}$ onto $\alpha < \cP(Y)$ let $A_f =
\{(y_1,y_2): f(y_1) < f(y_2)\} \subseteq Y \times Y$.
   Define $F:\cP(Y \times Y) \rightarrow \hrtg(Y)$ by $F(A) =
   \alpha$ if $A=A_f$ and $f,\alpha$ are as above, and $F(A) = 0$
 otherwise.

So $|\cP(\cP(Y)) \cup \bigcup\{{}^Y \alpha:\alpha < 
\hrtg([Y]^{\aleph_0})\}| \le_{\text{qu}} \cP(\cP(Y)) \times \cP(\cP(Y
\times Y))) = |\cP(\cP(Y \times Y))|$.  By the proof above we easily get
$|\text{Fil}^4_{\aleph_1}(Y)| \le_{\text{qu}} \cP(\cP(Y \times Y))$.
\end{PROOF}

\begin{claim}
\label{6.1}
[$\DC$]  Assume
\mn
\begin{enumerate}
\item[$(a)$]   ${\ga}$ is a countable set of limit ordinals
\sn
\item[$(b)$]   $<_*$ is a well ordering of $\Pi{\ga}$
\sn
\item[$(c)$]   $\theta \in {\ga} \Rightarrow { \text{\rm cf\/}}(\theta)
\ge \kappa$ where $\kappa = \hrtg({\cP}(\omega))$ or 
just $\Pi{\ga}/[{\ga}]^{< \aleph_0}$ is $< \kappa$-directed.
\end{enumerate}
\mn
\Then \, we can define $(\bar J,\bar{\gb},\bar{\bold f})$ such that
\mn
\begin{enumerate}
\item[$(\alpha)$] 
\begin{enumerate}
\item[(i)]  $\bar J = \langle J_i:i \le i(*) \rangle$ where $i(*) <
  \hrtg(\cP(\omega))$ 
\sn
\item[(ii)]  $J_i$ is an ideal on ${\ga}$ (though not necessarily a
  proper ideal)
\sn
\item[(iii)]  $J_i$ is increasing continuous with
$i,J_0 = \{\emptyset\},J_{i(*)} = {\cP}({\ga})$
\sn
\item[(iv)]  $\bar{\gb} = \langle {\gb}_i:i <
i(*) \rangle,{\gb}_i \subseteq {\ga}$ and $J_{i+1} = J_i + {\gb}_i$
\sn
\item[(v)]   so $J_i$ is the ideal on $\ga$ generated by $\{\gb_j:j <
  i\}$
\end{enumerate}
\sn
\item[$(\beta)$]
\begin{enumerate}
\item[(i)]  $\bar{\bold f} = \langle \bar f^i:i < i(*)\rangle$
\sn
\item[(ii)]  $\bar f^i = \langle f^i_\alpha:\alpha < \alpha_i \rangle$
\sn
\item[(iii)]  $f^i_\alpha \in \prod{\ga}$ is
$<_{J_i}$-increasing with $\alpha < \alpha_i$
\sn
\item[(iv)]  $\{f^i_\alpha:\alpha < \alpha_i\}$ is
cofinal in $(\prod{\ga},<_{J_i+({\ga} \backslash {\gb}_i)})$
\end{enumerate}
\sn
\item[$(\gamma)$]
\begin{enumerate}
\item[(i)]  $\cf(\prod{\ga}) \le \sum\limits_{i <i(*)} \alpha_i$
\sn
\item[(ii)]  for every $f \in \Pi{\ga}$ for some
$n$ and finite set $\{(i_\ell,\gamma_\ell):\ell < n\}$ such 
that $i_\ell < i(*),\gamma_\ell < \alpha_{i_\ell}$ we have 
$f < \max_{\ell < n}
f^{i_\ell}_{\gamma_\ell}$, i.e., $(\forall \theta \in 
{\ga})(\exists \ell < n)[f(\theta) < f^{i_\ell}_{\gamma_\ell}(\theta)]$.
\end{enumerate}
\end{enumerate}
\end{claim}

\begin{remark}
Note that there is no harm in having 
more than one occurence of $\theta \in \ga$.  See more 
in \cite{Sh:1005}, e.g. on uncountable $\ga$.
\end{remark}

\begin{PROOF}{\ref{6.1}}
Note that:
\mn
\begin{enumerate}
\item[$\circledast_1$]   clause $(\gamma)$ follows from $(\alpha) +
(\beta)$.
\end{enumerate}
\mn
[Why?  Easily $(\gamma)(ii) \Rightarrow (\gamma)(i)$.  Now let $g \in
\Pi{\ga}$ and let $I_g = \{\gb \subseteq \ga$: we can find $n < \omega$
and $i_\ell < i(*)$ and $\beta_\ell < \alpha_{i_\ell}$ for $\ell < n$
such that $\theta \in {\gb} \Rightarrow (\exists \ell <
n)(g(\theta) < f^{i_\ell}_{\beta_\ell}(\theta))\}$.

Easily $I_g$ is an ideal on $\ga$ though not necessarily a proper
ideal.  Note that if ${\ga} \in
I_g$ we are done.  So assume ${\ga} \notin I_g$.  Note that 
$I_g \subseteq J_{i(*)}$ hence $j_g = \text{ min}\{i \le i(*)$: some
${\gc} \in {\cP}({\ga}) \backslash I_g$ belongs to $J_i\}$
is well defined (as ${\ga} \in {\cP}({\ga}) \backslash I_g
\wedge {\ga} \in J_{i(*)}$).  As $J_0 = \{\emptyset\}$ and clearly
if $\emptyset \in \cT_g$ we have $j_a >0$.  As
$\langle J_i:i \le i(*)\rangle$ is $\subseteq$-increasing continuous,
necessarily $j_g$ is a successor ordinal say $j_g = i_g +1$ and let
$i(g) = i_g$ and choose $\gc \in J_{i(g)} \backslash I_g$, clearly
$J_{i(g)} \subseteq I_g$ so $\gc$ belongs to
$J_{j_g} \backslash J_{i_g}$.  
By clause $(\beta)(iv)$ there is $\alpha < \alpha_{i(g)}$
such that $g < f^i_\alpha \mod(J_{i(g)} + (\ga \backslash
\gb_{i(g)}))$.

Now let $\gd = \{\theta \in \ga:g(\theta) <
f^i_\alpha (\theta)\}$ so by the choice of $\alpha$ we have
$\gd = \ga \mod (J_{i(g)} + (\ga \backslash \gb_{(g)})$ which means
that $\gb_{i(g)} \subseteq \gd \mod J_{i(g)}$ so as $J_{i(g)+1} =
J_{i(g)} + \gb_{i,g}$ and $\gc \in J_{i(g)+1} \backslash J_{i(g)}$
clearly $\gc \subseteq \gb_{i(g)} \mod J_{i(g)}$.

But by the definition of the ideal $J_{i(g)}$ and of $\gd$ necessarily $\gd
\in J_{i(g)}$ and recall $J_{i(g)} \subseteq J_{i(g)}$, contradicting the
conclusion of the last sentence.]

Since $(\gamma)$ follows from $(\alpha) + (\beta)$, it suffices to
prove these parts. 
By induction on $i < \kappa$ we try to choose $(\bar J^i,\bar{\gb}^i,
\bar{\bold f}^i)$ where $\bar J^i = \langle J_j:j
\le i \rangle,\bar{\gb}^i = \langle {\gb}^i_j:
j < i\rangle,\bar{\bold f}^i = \langle \bar f^j:
j < i \rangle$ which satisfies the relevant
parts of the conclusion and do it uniformly from $({\ga},<_*)$.  
Once we arrive at $i$ such that $J_i = {\cP}({\ga})$ we are done.

For $i=0$ recalling $J_0 = \{\emptyset\}$ there is no problem.

For $i$ limit recalling that $J_i = \cup\{J_j:j<i\}$ there is no
problem and note that if $j<i \Rightarrow {\ga} \notin J_j$ then
${\ga} \notin J_i$.

So assume that $(\bar J^i,\gb^i,\bar{\bold f}^i)$ is well defined and
$\ga \notin J_i$ and we shall define for $i+1$.

We try to choose $\bar g^{i,\varepsilon} = \langle
g^{i,\varepsilon}_\alpha:\alpha < \delta_{i,\varepsilon}\rangle$ and
$\gb_{i,\varepsilon}$ by
induction on $\varepsilon < \omega_1$ and for each $\varepsilon$ we
try to choose $g^{i,\varepsilon}_\alpha \in \Pi{\ga}$ by induction
on $\alpha$ (in fact $\alpha <$ \hrtg$(\Pi{\ga})$ suffice, we
shall get stuck earlier) such that:
\mn
\begin{enumerate}
\item[$\circledast^2_{i,\varepsilon}$]
\begin{enumerate}
\item[(a)]   if $\beta < \alpha$ then $g^{i,\varepsilon}_\beta <_{J_i}
g^{i,\varepsilon}_\alpha$
\sn
\item[(b)]  if $\zeta < \varepsilon$ and $\alpha <
\delta_{i,\zeta}$ then $g^{i,\zeta}_\alpha \le g^{i,\varepsilon}_\alpha$
\sn
\item[(c)]  if cf$(\alpha) = \aleph_1$ then
$g^{i,\varepsilon}_\alpha$ is defined by 
\[
\theta \in {\ga} \Rightarrow g^{i,\varepsilon}_\alpha(\theta) =
\text{ Min}\{\bigcup\limits_{\beta \in C} g^{i,\varepsilon}_\beta(\theta):C
\text{ is a club of } \alpha\}
\]
\sn
\item[(d)]  if $\alpha$ is a limit ordinal and
cf$(\alpha) \ne \aleph_1,\alpha \ne
0$ then $g^{i,\varepsilon}_\alpha$ is the $<_*$-first $g \in 
\Pi{\ga}$ satisfying clauses (a) + (b)
\sn
\item[(e)] if we have $\langle
g^{i,\varepsilon}_\beta:\beta < \alpha \rangle$, cf$(\alpha) >
\aleph_1$, moreover $\cf(\alpha) \ge \min\{\cf(\theta):\theta \in \ga\}$ 
and there is no $g$ as required in clause (d) then 
$\delta_{i,\varepsilon} = \alpha$
\sn
\item[(f)]   if $\alpha=0$ or $\alpha$ is a successor, 
then $g^{i,\varepsilon}_\alpha$ is the $<_*$-first $g \in \Pi{\ga}$
such that:
\sn
\begin{enumerate}
\item[$\bullet_1$]  $\zeta < \varepsilon \wedge \alpha < \delta_{i,\zeta} 
\Rightarrow g^{i,\zeta}_\alpha \le g$ 
\sn
\item[$\bullet_2$]  $\beta < \alpha \Rightarrow g^{i,\varp}_\beta <
  g^{i,\varp}_\alpha \mod J_i$
\sn
\item[$\bullet_3$]  $\varepsilon = \zeta +1 \Rightarrow (\forall \beta <
\delta_{i,\zeta})[\neg(g \le_{J_i} g^{i,\zeta}_\beta)]$, follows if 
$\alpha > 0$
\end{enumerate}
\sn
\item[(g)]  $J_i$ is the ideal on $\cP(\ga)$
 generated by $\{\gb_j:j < i\}$
\sn
\item[(h)]  $\gb_{i,\varepsilon} \in (J_i)^+$ so
 $\gb_{i,\varepsilon} \subseteq \ga$
\sn
\item[(i)]  $\bar g^{i,\varepsilon}$ is increasing and
 cofinal in $(\Pi(\ga),<_{J_i + (\ga \backslash \gb_{i,\varepsilon})})$
\sn
\item[(j)]  $\gb_{i,\varepsilon}$ is such that under
clauses $(h) + (i)$ the set $\{\text{otp}(\ga \cap \theta):\theta \in
\gb_{i,\varepsilon}\}$ is $<_*$-minimal
\sn
\item[(k)]  $\gb_{i,\zeta} \subseteq \gb_{i,\varepsilon}$ mod $J_i$
  (follows by ``if $\zeta < i$ then $g^{i,\varepsilon}_0$
 is a $<_{J_i + \gb_{i,\zeta}}$-upper bound of $\bar g^{i,\zeta}$".
\end{enumerate}
\end{enumerate}
\mn
Clearly in stage $\varp$ we first choose $g^{i,\varp}_\alpha$ by
induction on $\alpha$.  As $\beta < \alpha \Rightarrow
g^{i,\varp}_\beta \ne g^{i,\varp}_\alpha$ we are stuck in some
$\delta_{i,\varp}$ and then choose $\gb_{i,\varp}$.

We now give details on some points:
\mn
\begin{enumerate}
\item[$(*)_0$]  if $\alpha = 0$ then we can choose $g^{2,\varp}_0$.
\end{enumerate}
\mn
[Why?  Trivial.]
\mn
\begin{enumerate}
\item[$(*)_1$]  Clause (c) is O.K., that is:  if we arrive to
  $(\varp,\alpha),\cf(\alpha) = \aleph_1$ then we can define
  $g^{i,\varp}_\alpha$.
\end{enumerate}
\mn
[Why?  We already have $\langle g^{i,\varepsilon}_\alpha:\alpha <
\delta\rangle$ and $\langle g^{i,\zeta}_\alpha:\alpha <
\delta_{i,\zeta},\zeta < \varepsilon\rangle$, and we define
$g^{i,\varepsilon}_\delta$ as there.  Now
$g^{i,\varepsilon}_\delta(\theta)$ is well defined as the ``Min" is
taken on a non-empty set of ordinals  
as we are assuming cf$(\delta) = \aleph_1$.
The value is $< \theta$ because for some club $C$ of $\delta$, otp$(C)
= \omega_1$, so $g^{i,\varepsilon}_\delta(\theta) \le
\cup\{g^{i,\varepsilon}_\beta(\theta):\beta \in C\}$ but this set is
$\subseteq \theta$ while cf$(\theta) > \aleph_1$ by clause (c) of the 
assumption.
By AC$_{\aleph_0}$ we can find a sequence $\langle C_\theta:\theta \in
{\ga}\rangle$ such that: $C_\theta$ is a club of $\delta$ of order
type $\omega_1$ satisfying $g^{i,\varepsilon}_\delta(\theta) =
\cup\{g^{i,\varepsilon}_\alpha(\theta):\alpha \in C_\theta\}$ hence
for every club $C$ of $\delta$ included in $C_\theta$ we have
$g^{i,\varepsilon}_\delta(\theta) =
\cup\{g^{i,\varepsilon}_\alpha(\theta):\alpha \in C_\theta\}$.  Now
$\theta \in {\ga} \Rightarrow g^{i,\varepsilon}_\delta(\theta) 
= \bigcup\limits_{\alpha \in C}
g^{i,\varepsilon}_\alpha(\theta)$ when $C := 
\cap\{C_\sigma:\sigma \in {\ga}\}$, because $C$ too is a club of
$\delta$ recalling ${\ga}$ is countable.  So if $\alpha < \delta$
then for some $\beta$ we have 
$\alpha < \beta \in C$ hence the set ${\gc} := \{\theta
\in {\ga}:g^{i,\varepsilon}_\alpha(\theta) \ge
g^{i,\varepsilon}_\beta(\theta)\}$ belongs to $J_i$ and $\theta \in
\ga \backslash \gc \Rightarrow
g^{i,\varepsilon}_\alpha(\theta) < g^{i,\varepsilon}_\beta(\theta) \le
g^{i,\varepsilon}_\delta(\theta)$, so indeed $g^{i,\varepsilon}_\alpha
<_{J_i} g^{i,\varepsilon}_\delta$. 

Lastly, why $\zeta < \varepsilon \Rightarrow g^{i,\zeta}_\delta \le
g^{i,\varepsilon}_\delta$?  As we can find a club $C$ of $\delta$
which is as above for both $g^{i,\zeta}_\delta$ and
$g^{i,\varepsilon}_\delta$ and recall that clause (b) of
$\circledast_{i,\varepsilon}$ holds for every $\beta \in C$.   Together 
$g^{i,\varepsilon}_\delta$ is as required.]
\mn
\begin{enumerate}
\item[$(*)_2$]    cf$(\delta_{i,\varp}) > \aleph_1$ and even
cf$(\delta_{i,\varp}) \ge \text{ min}\{\text{cf}(\theta):\theta \in \ga\}$.
\end{enumerate}
\mn
[Why?  We have to prove that arriving to $\alpha > 0$, if cf$(\alpha) <
\text{ min}\{\text{cf}(\theta):\theta \in \ga\}$ then we can choose
$g^{i,\varp}_\alpha$ as required.  The cases cf$(\alpha) =
\aleph_1,\alpha=0$ are covered by $(*)_1,(*)_0$ respectively, 
otherwise let $u \subseteq \alpha$ be
unbounded of order type cf$(\alpha)$, and define a function $g$ from
$\ga$ to the ordinals by $g(\theta) = \sup(\{g^{i,\varp}_\beta(\theta):\beta
\in u\} \cup \{g^{i,\zeta}_\alpha(\theta):\zeta < \varepsilon\})$.
This is a subset of $\theta$ of cardinality $< |\ga| + \text{
cf}(\alpha)$ which is $< \theta = \text{ cf}(\theta)$ hence
$g \in \Pi\ga$, easily is as required, i.e. satisfies clauses (a) + (b)
and the $<_*$-first such $g$ is $g^{i,\varp}_\alpha$.]

Note that clause (e) of $\circledast_{i,\varepsilon}$ follows.
\mn
\begin{enumerate}
\item[$(*)_3$]   if $\zeta < \varepsilon$ then $\delta_{i,\varepsilon}
\le \delta_{i,\zeta}$.
\end{enumerate}
\mn
[Why?  Otherwise $g^{i,\varepsilon}_{\delta_{i,\zeta}}$ contradict
clause (e) of $\circledast_{i,\zeta}$.]
\mn
\begin{enumerate}
\item[$(*)_4$]   if $g^{i,\varp} = \langle g^{i,\varp}_\alpha:\alpha <
  \delta_{i,\varp}\rangle$ is well defined and $\cf(\delta_{i,\varp})
  ge \kappa$ then $\gb_{i,\varp}$ is well defined.
\end{enumerate}
\mn
[Why?  Clearly it suffices to prove that there is $\gb$ as required on
$\gb_{i,\varp}$ (in clauses (b),(i)).  So toward contradiction assume
that for every $\gb \in J^+_i,\bar g^{i,\varp}$ is not
$<_{J_i}$-cofinal in $\Pi \ga$ hence there is $h \in \Pi \ga$ such
that $\alpha < \delta_{i,\varp} \Rightarrow h \nleq_{J_i}
g^{i,\varp}_\alpha$ and let $h_b$ be the $<_*$-minimal such $h$.  Let
$h_*$ be the function with domain $\ga$ such that $h(\theta)
=\cup\{h_{\gb}(\theta) +1:b \in J^+_i\}$.

As $\hrtg(J^+_i) \le \hrtg(\cP(\ga)) < \min\{\cf(\theta):\theta \in
\ga\}$, clearly $h_* \in \Pi \ga$.  Now for $\alpha <
\delta_{i,\varp}$ let $\gd_{i,\varp,\alpha} = 
\{\theta \in \ga:g^{i,\varp}_\alpha(\theta) \le h_*(\theta)\}$.  So
$\langle \gd_{i,\varp,\alpha}/J_i:\alpha < \delta_{i,\varp}\rangle$ is
$\le$-increasing in the Boolean Algebra $\cP(\ga)/J_i$, so for some
$\beta_{i,\varp} < \delta_{i,\varp}$ we have $\alpha \in
(\beta_{i,\varp},\delta_{i,\varp}) \Rightarrow \gd_{i,\varp,\alpha} =
\gd_{i,\varp,\beta_{i,\varp}} \mod J_i$.  This implies
$\gd_{i,\varp}$ can serve as $\gb_{i,\varp}$.]

To finish consider the following two cases.
\bigskip

\noindent
\underline{Case 1}:  We succeed to carry the induction, i.e. choose
$\bar g^{i,\varepsilon}$ for every $\varepsilon < \kappa$.

So $\langle \gb_{i,\varp}:\varp < \kappa\rangle$ is a sequence of
subsets of $\ga$, pairwise distinct (by $\circledast^2_{\kappa,0}$
clauses (g) + (b)), but $\kappa \ge \hrtg(\cP(\omega))$ and $\ga$ is
countable; contradiction.
\bigskip

\noindent
\underline{Case 2}:  We are stuck in $\varepsilon < \kappa$.

For $\varepsilon = 0$ there is no problem to define
$g^{i,\varepsilon}_\alpha$ by induction on $\alpha$ till we are stuck,
say in $\alpha$, necessarily $\alpha$ is of large enough cofinality
$\ge \kappa$ by $(*)_2$, and so 
$\bar g^{i,\varepsilon}$ is well defined.  We then prove
$\gb_{i,\varepsilon}$ exists by $(*)_4$ again using $<_*$.

For $\varepsilon$ limit we can also choose $\bar g^\varepsilon$.

For $\varepsilon = \zeta +1$, if $\ga \in J_\varp$ then we are done;
otherwise $g^{i,\varepsilon}_0$ as required can be chosen by $(*)_0$, and then
we can prove that 
$\bar g^{i,\varepsilon},\gb_{i,\varepsilon}$ exists as above. 
\end{PROOF}

\begin{remark}
\label{6.2p}
\relax From \ref{6.1} we can deduce bounds on
$\hrtg({}^Y(\aleph_\delta))$ when $\delta < \aleph_1$ and more like
the one on $\aleph^{\aleph_0}_\omega$ (better the bound on
pp$(\aleph_\omega))$. 
\end{remark}
\newpage

\section {No decreasing sequence of subalgebras} \label{2}

In this section we 
concentrate on weaker axioms.  We consider Theorem \ref{r.2}
under weaker assumptions than ``$[\lambda]^{\aleph_0}$ is well
orderable".  We are also interested in replacing $\omega$ by
$\partial$ in ``no decreasing $\omega$-sequence of $c \ell$-closed
sets", but the reader may consider $\partial = \aleph_0$ only.  Note
that for the full version, Ax$^4_\alpha$, i.e., $[\alpha]^\partial$ is
well orderable, the case of $\partial = \aleph_0$ is implied
by the $\partial > \aleph_0$ version and suffices for the results.
But for other versions, the axioms for different $\partial$'s seem
incomparable. 

Note that if we add many Cohens (not well ordering them) then
Ax$^4_\lambda$ fails below even for $\partial =  \aleph_0$, whereas
the other axioms are not affected.  But forcing by $\aleph_1$-complete
forcing notions preserve Ax$_4$.

\begin{hypothesis}
\label{22.0}
DC$_\partial$ and let $\partial(*) = \partial + \aleph_1$.  Actually
we use only DC in \ref{r.4}(1) and DC$_\partial$ in \ref{r.4}(3) and
the later claims.  We fix a regular cardinal $\partial$.
\end{hypothesis}

\begin{definition}
\label{22.1}
Below we should, e.g. write Ax$^{\ell,\partial}$ instead of 
Ax$^\ell$ and assume $\alpha > \mu > \kappa \ge \partial$.
If $\kappa = \partial$ we may omit it. 

\noindent
1) Ax$^1_{\alpha,\mu,\kappa}$ means that
there is a weak closure operation on $\lambda$ of character
$(\mu,\kappa)$, see Definition \ref{0.F}(1A),
such that there is no $\subseteq$-decreasing $\partial$-sequence
$\langle \cU_\varepsilon:\varepsilon <\partial\rangle$ of
subsets of $\alpha$ with $\varepsilon <\partial \Rightarrow c
\ell(\cU_{\varepsilon +1}) \nsupseteq \cU_\varepsilon$.  
We may here and below replace
$\kappa$ by $<\kappa$; similarly for $\mu$; let $< |Y|^+$ means $|Y|$.

\noindent
2) Let Ax$^0_{\alpha,<\mu,\kappa}$ mean there is a function 
$c \ell:[\alpha]^{\le \kappa} \rightarrow [\alpha]^{< \mu}$
such that $u \cup \{0\} \subseteq c \ell(u)$ and there is no
$\subseteq$-decreasing sequence $\langle {\cU}_\varepsilon:
\varepsilon < \partial \rangle$ of members of $[\alpha]^{\le \kappa}$ 
such that $\varepsilon < \partial
\Rightarrow c \ell({\cU}_{\varepsilon +1}) \nsupseteq {\cU}_\varepsilon$. 

\noindent
2A) Writing $Y$ instead of $\kappa$ means 
$c \ell:[\alpha]^{< \hrtg(Y)} \rightarrow [\alpha]^{<\mu}$.  Let $c
\ell_{[\varepsilon]}:{\cP}(\alpha) \rightarrow {\cP}(\alpha)$ be
$c \ell^1_{\varepsilon,< \text{reg}(\kappa^+)}$ as defined in
\ref{0.F}(4) recalling reg$(\gamma) = \text{ Min}\{\chi:\chi$ a regular
cardinal $\ge \gamma\}$.

\noindent
3) Ax$^2_\alpha$ means that there is ${\cA} \subseteq
[\alpha]^\partial$ which is well orderable and for every $u \in
[\alpha]^\partial$ for some $v \in {\cA},u \cap v$ has power $= \partial$. 

\noindent
4) Ax$^3_\alpha$ means that cf$([\alpha]^{\le \partial},\subseteq)$ is
below some cardinal, i.e., some cofinal ${\cA} \subseteq
[\alpha]^\partial$ (under $\subseteq$) is well orderable.

\noindent
5) Ax$^4_\alpha$ means that $[\alpha]^{\le \partial}$ is well orderable.

\noindent
6) Above omitting $\alpha$ (or writing $\infty$) means ``for every
$\alpha$", omitting $\mu$ we mean ``$< \hrtg({\cP}(\partial))$".

\noindent
7) Lastly, let Ax$_\ell = \text{ Ax}^\ell$ for $\ell = 1,2,3$.
\end{definition}

\noindent
So easily (or we have shown in the proof of \ref{r.2}):
\begin{claim}
\label{r.2b}
1) {\rm Ax}$^4_\alpha$ implies {\rm Ax}$^3_\alpha$, 
{\rm Ax}$^3_\alpha$ implies {\rm Ax}$^2_\alpha$, 
{\rm Ax}$^2_\alpha$ implies {\rm Ax}$^1_\alpha$ and Ax$^1_\alpha$ implies
Ax$^0_\alpha$ .  Similarly for {\rm Ax}$^\ell_{\alpha,< \mu,\kappa}$.

\noindent
2) In Definition \ref{22.1}(2), the last demand,
if $c \ell$ has monotonicity, \then \, only
$c \ell \restriction [\alpha]^{\le \partial}$ is relevant, in fact, an
equivalent demand is that if $\langle \beta_\varepsilon:\varepsilon <
\partial \rangle \in{}^\partial \alpha$ then for some
$\varepsilon,\beta_\varepsilon \in c \ell\{\beta_\zeta:\zeta \in
(\varepsilon,\partial)\}$.

\noindent
3) If {\rm Ax}$^0_{\alpha,< \mu_1,< \theta}$ and $\theta \le$ \hrtg$(Y)$ and
\footnote{Can do somewhat better; we can replace $[\alpha]^{< \mu_1}$
  by $\{v \subseteq \alpha:\otp(v) \subseteq \mu_1\}$} 
$\mu_2 = \sup\{\hrtg(\mu_1 \times
[\beta]^\theta):\beta < \hrtg(Y)\}$ \then \, 
Ax$^0_{\alpha,< \mu_2,< \hrtg(Y)}$. 
\end{claim}

\begin{PROOF}{\ref{r.2b}}
1) Clearly Ax$^2_{\alpha,< \mu,\kappa} \Rightarrow
\text{ Ax}^1_{\alpha,< \mu,\kappa}$ holds similarly to the proof of
\ref{r.2y}; the other implications hold by inspection.

\noindent
2) First assume that we have a $\subseteq$-decreasing sequence
$\langle \cU_\varepsilon:\varepsilon < \partial\rangle$ such that
   $\varepsilon < \partial \Rightarrow c \ell(\cU_{\varepsilon +1})
\nsupseteq \cU_\varepsilon$.  Let $\beta_\varepsilon = \text{
 min}(\cU_\varepsilon \backslash c \ell(\cU_{\varepsilon +1}))$ for
   $\varepsilon < \partial$ so clearly $\bar\beta = \langle
\beta_\varepsilon:\varepsilon < \partial \rangle$ exists; so by
monotonicity $c \ell(\{\beta_\zeta:\zeta \in [\varepsilon
+1,\partial)\} \subseteq c \ell(\cU_{\varepsilon +1})$ hence
$\beta_\varepsilon \notin c \ell(\{\beta_\zeta:\zeta \in
   [\varepsilon +1,\partial)\}$.

Second, assume that $\bar\beta = \langle \beta_\varepsilon:\varepsilon
< \partial\rangle \in {}^\partial \alpha$ satisfies $\beta_\varepsilon
\notin c \ell(\{\beta_\zeta:\zeta \in [\varepsilon +1,\partial)\}$ for
$\varepsilon < \partial$.  Now letting $\cU'_\varepsilon =
\{\beta_\zeta:\zeta < \partial$ satisfies 
$\varepsilon \le \zeta\}$ for $\varepsilon < \partial$
clearly $\langle \cU'_\varepsilon:\varepsilon < \partial\rangle$ exists, is
$\subseteq$-decreasing and $\varepsilon < \partial \Rightarrow
\beta_\varepsilon \notin c \ell(\cU'_{\varepsilon +1}) \wedge
\beta_\varepsilon \in \cU'_\varepsilon$.  So we have
shown the equivalence.

\noindent
3) Let $c \ell(-)$ witness Ax$^0_{\alpha,< \mu_1,< \theta}$.
We define the function $c \ell'$ with domain 
$[\alpha]^{< \hrtg(Y)}$ by $c \ell'(u) = \cup\{c \ell(v):v \subseteq u$ has
cardinality $< \theta\}$.

Now
\mn
\begin{enumerate}
\item[$(*)_0$]   $c \ell'$ is a function from $[\alpha]^{< \hrtg(Y)}$
  into $[\alpha]^{< mu_2}$.
\end{enumerate}
\mn
For this it is enough to note:
\mn
\begin{enumerate}
\item[$(*)_1$]   if $u \in [\alpha]^{< \hrtg(Y)}$ then
$c \ell'(u)$ has cardinality $< \mu_2 := \sup\{\hrtg(\mu_1 \times
[\beta]^\theta:\beta < \hrtg(Y)\}$.
\end{enumerate}
\mn
[Why?  Let $C_u = \{(v,\varepsilon):v \subseteq u$ has cardinality $<
\theta$ and $\varepsilon < \text{ otp}(c \ell(v))$ which is $<
\mu_1\}$.  Clearly $|c \ell'(u)| < \hrtg(C_u)$ and $|C_u| = |\mu_1
\times [\text{otp}(u)]^{<\theta}|$, so $(*)_1$ holds.  Note that if
$\alpha_* < \mu^+_1$ we can replace the demand $v \in [u]^{< \theta}
\Rightarrow |c \ell(v)| < \mu_1$ by $v \in [u]^{< \theta} \Rightarrow
\text{ otp}(c \ell(v)) < \alpha_*$.] 
\mn
\begin{enumerate}
\item[$(*)_2$]   If $\langle u_\varepsilon:\varepsilon < \partial
\rangle$ is $\subseteq$-decreasing where $u_\varp \subseteq \alpha$
then $u_\varepsilon \subseteq c
\ell'(u_{\varepsilon +1})$ for some $\varepsilon < \partial$.
\end{enumerate}
\mn
[Why?  If not we can choose a sequence $\langle
\beta_\varepsilon:\varepsilon < \partial \rangle$ by letting
$\varepsilon <\partial \Rightarrow \beta_\varepsilon = \text{ min}
(u_\varepsilon \backslash c \ell'(u_{\varepsilon +1}))$.  Let
$u'_\varepsilon = \{\beta_\zeta:\zeta \in [\varepsilon,\partial)\}$.
As $\langle u'_\varp:\varp < \partial\rangle$ is $\subseteq$-decreasing
by the choice of $c \ell(-)$ for some
$\varepsilon,\beta_\varepsilon \in c \ell\{\beta_\zeta:\zeta \in
(\varp +1,\partial)\}$, but this set is $\subseteq c
\ell'(u_{\varepsilon +1})$ by the definition of $c \ell'(-)$, 
so we are done.]  
\end{PROOF}

\begin{claim}
\label{r.B.3}
Assume $c \ell$ witness 
{\rm Ax}$^0_{\alpha,< \mu,\kappa}$ so $\partial \le \kappa < \mu$ and so
$c \ell:[\alpha]^{\le \kappa} \rightarrow [\alpha]^{< \mu}$ 
and recall $c \ell^1_{\varepsilon,\le \kappa}:\cP(\alpha) 
\rightarrow {\cP}(\alpha)$ is from \ref{22.1}(2A), \ref{0.F}(4).

\noindent
1) $c \ell^1_{1,\le \kappa}$ is a weak closure operation, it has character
$(\mu_\kappa,\kappa)$ whenever $\partial \le \kappa \le \alpha$ 
and $\mu_\kappa =$ \hrtg$(\mu \times {\cP}(\kappa))$, see
Definition \ref{0.F}.

\noindent
2) $c \ell^1_{\text{\rm reg}(\kappa^+),\le \kappa}$ is a 
closure operation and it has character $(< \mu'_\kappa,\kappa)$ when $\partial
\le \kappa \le \alpha$ and $\mu'_\kappa =$ 
\hrtg$({\cH}_{<\partial^+}(\mu \times \kappa))$.
\end{claim}

\begin{PROOF}{\ref{r.B.3}}
1) By its definition $c \ell^1_{1,\le \kappa}$ is a weak closure operation.

Assume $u \subseteq \alpha,|u| \le \kappa$; non-empty for simplicity.
Clearly $\mu \times [|u|]^{< \partial}$ has the same power as $\mu \times
[u]^{< \partial}$.  Define \footnote{clearly we can replace $< \mu$ by $<
\gamma$ for $\gamma \in (\mu,\mu^+)$}
 the function $G$ with domain $\mu \times [u]^{< \partial}$ as 
follows: if $\alpha < \mu$ and $v \in [u]^{\le
\partial}$ then $G((\alpha,v))$ is the $\alpha$-th member of $c \ell(v)$
if $\alpha < \text{ otp}(c \ell(v))$ and $G((\alpha,v)) = \text{
min}(u)$ otherwise.

So $G$ is a function from $\mu \times [u]^{\le \partial}$ onto $c
\ell^1_{1,\le \kappa}(u)$.  This proves that $c \ell^1_{1,\le \kappa}$
has character $(< \mu_\kappa,\kappa)$ as $\mu_\kappa =$ 
\hrtg$(\mu \times {\cP}(\kappa))$.

\noindent
2) If $\langle u_\varepsilon:\varepsilon \le \text{ reg}(\kappa^+)\rangle$ is
an increasing continuous sequence of sets \then \, $[u_{\partial^+}]^{\le
\partial} = \cup\{[u_\varepsilon]^{\le \partial}:\varepsilon <
\text{ reg}(\kappa^+)\}$ as $\text{\rm reg}(\kappa^+)$ is regular
(even of cofinality $> \partial$ suffice)
by its definition, note reg$(\partial^+) = \partial^+$ when
AC$_\partial$ holds when DC$_\partial$ holds.

Second, let $u \subseteq \alpha,|u| \le \kappa$ and let $u_\varepsilon
= c \ell^1_{\varepsilon,\kappa}(u)$ for $\varepsilon \le \partial^+$; it is
enough to show that $|u_{\partial^+}| < \mu'_\kappa$.  The proof is
similar to earlier one.  
\end{PROOF}

\begin{dc}
\label{r.4}  
Let $c \ell$ exemplify Ax$^0_{\lambda,< \mu,Y}$ and 
$Y$ be an uncountable set such that $\partial(*) \le_{\qu} Y$.

\noindent
1) Let ${\cF}_{\gy},{\cF}_{{\gy},\alpha}$ be as in the proof
of Theorem \ref{r.2} for ${\gy} \in \text{ Fil}^4_{\partial(*)}(Y,\mu)$ 
and ordinal $\alpha$ (they depend on $\lambda$ and $c \ell$ but note
that $c \ell$ determines $\lambda$; so if we
derive $c \ell$ by Ax$^4_\lambda$ then they depend indirectly on the well 
ordering of $[\lambda]^\partial$) so we may write $\cF_{\gy,\alpha}
 = {\cF}_{\gy}(\alpha,c \ell)$, etc.

That is, fully
\mn
\begin{enumerate}
\item[$(*)_1$]   for ${\gy} \in \text{ Fil}^4_{\partial(*)}(Y,\mu)$
and ordinal $\alpha$ let ${\cF}_{{\gy},\alpha}$ be the set of $f$ such that:
\sn
\begin{enumerate}
\item[$(a)$]   $f$ is a function from $Z^{\gy}$ to $\lambda$
\sn
\item[$(b)$]   $\text{\rm rk}_{D[\gy]}(f) = \alpha$ recalling
that this means rk$_{D^{\gy}_1 + Z^{\gy}}(f \cup 0_{Y \backslash
  Z^{\gy}}) = \alpha$ by Definition \ref{0.A}(2)
\sn
\item[$(c)$]    $D^{\gy}_2 = D^{\gy}_1 \cup \{Y \backslash
A:A \in J[f,D^{\gy}_1]\}$, see Definition \ref{0.C}
\sn
\item[$(d)$]   $Z^{\gy} \in D^{\gy}_2$
\sn
\item[$(e)$]   if $Z \in D^{\gy}_2$ and $Z \subseteq 
Z^{\gy}$ then $c \ell(\{f(y):y \in Z\}) \supseteq \{f(y):y \in Z^{\gy}\}$
\sn
\item[$(f)$]   $h^{\gy}$ is a function with domain $Z^{\gy}$ 
such that $y \in Z^{\gd} \Rightarrow h^{\gh}(y) = \text{ otp}
(f(y) \cap \{c \ell(\{f(z):z \in Z^{\gy}\})$
\end{enumerate}
\item[$(*)_2$]   ${\cF}_{\gy} = \cup\{{\cF}_{{\gy},\alpha}:
\alpha$ an ordinal$\}$.
\end{enumerate}
\mn
2) Notice that $\cF_{\gy,\alpha}$ is a singleton or the empty set.
Let $\Xi_{\gy} = \Xi_{\gy}(c \ell) = \Xi_{\gy}(\lambda,c \ell)
 = \{\alpha:{\cF}_{{\frak y},\alpha} \ne \emptyset\}$ and 
$f_{\gy,\alpha}$ is the function 
$f \in {\cF}_{{\gy},\alpha}$ when $\alpha \in \Xi_{\gy}$; it is well defined.

\noindent
3) If $D \in \text{ Fil}_{\partial(*)}(Y)$, rk$_D(f) = \alpha$ and $f
\in {}^Y \lambda$ \then \, $\alpha \in \Xi_D(\lambda,c \ell)$
and $f \restriction Z^{\gy} = f_{{\gy},\alpha}$ for some $\gy \in
\text{ Fil}^4_{\aleph_1}(Y)$; moreover, $(D^{\gy}_1,D^{\gy}_2) =
(D,\dual(J(J[f,D]))$ where $\Xi_D(\lambda,c \ell)
:= \cup\{\Xi_{\gy}:{\gy} \in \text{ Fil}^4_{\partial(*)}(Y)$ and
$D^{\gy}_1 = D\}$.

\noindent
4) If $D \in \text{ Fil}_{\partial(*)}(Y),
f \in {}^Y \lambda,Z \in  D^+$ and rk$_{D+Z}
(f) \ge \alpha$ \then \, for some $g \in
\prod\limits_{y \in Y} (f(y)+1) \subseteq {}^Y(\lambda +1)$ we have
$\text{rk}_D(g) = \alpha$ hence $\alpha \in \Xi_D(\lambda,c \ell)$.

\noindent
5) So we should write ${\cF}_{\gy}[c \ell],
\Xi_{\gy}[\lambda,c \ell],f_{{\gy},\alpha}[c \ell]$.
\end{dc}

\begin{PROOF}{\ref{r.4}}
As in the proof of \ref{r.2} recalling ``$c \ell$ exemplifies
Ax$^0_{\lambda,< \mu,\hrtg(Y)}$" holds, this replaces the use of $F_*$ there; 
and see the proof of \ref{r.9} below in part (3), for this we need:
\mn
\begin{enumerate}
\item[$\boxplus$]  if $D \in \Fil^1_\partial(Y)$ and $f \in
  {}^\kappa \partial$, \then \, for some $Z \in D$ we have:
\sn
\begin{itemize}
\item  if $Y \subseteq Z$ belongs to $D$ then $c \ell(\Rang(f \rest Y)
  = c \ell(\Rang(f \rest Z))$.
\end{itemize}
\end{enumerate}
\mn
[Why $\boxplus$ holds?  By Definition \ref{22.1}(2) using the axiom
DC$_\partial$.] 
\end{PROOF}

\begin{claim}
\label{y.21}
We have $\xi_2$ is an ordinal and {\rm Ax}$^0_{\xi_2,< \mu_2,Y}$ holds
\when ,(note that $\mu_2$ is not much larger than $\mu_1$):
\mn
\begin{enumerate}
\item[$(a)$]    Ax$^0_{\xi_1,< \mu_1,Y}$ so $\partial < \hrtg(Y)$
\sn
\item[$(b)$]   $c \ell$ witnesses clause (a)
\sn
\item[$(c)$]   $D \in { \text{\rm Fil\/}}_{\partial(*)}(Y)$
\sn
\item[$(d)$]   $\xi_2  = 
\{\alpha:f_{{\gy},\alpha}[c \ell]$ is well defined 
for some ${\gy} \in { \text{\rm Fil\/}}^4_{\partial(*)}(Y,\mu_1)$
which satisfies $D^{\gy}_1 = D$ and necessarily {\rm Rang}$(f_{{\gy},
\alpha}[c \ell]) \subseteq \xi_1\}$
\sn
\item[$(e)$]   $\mu_2$ is defined as $\mu_{2,3}$ where:
\sn
\begin{enumerate}
\item[$(\alpha)$]   let $\mu_{2,0} =$ \hrtg$(Y)$
\sn
\item[$(\beta)$]   $\mu_{2,1} = \sup_{\beta < \mu_{2,0}} \hrtg(\beta \times
\text{\rm Fil}^4_{\partial(*)}(Y,\mu_1))$
\sn
\item[$(\gamma)$]   $\mu_{2,2} = \sup_{\alpha < \mu_{2,1}}
\hrtg(\mu_1 \times [\alpha]^{\le\partial})$
\sn
\item[$(\delta)$]    $\mu_{2,3} = \sup\{\hrtg({}^Y \beta \times 
\Fil_{\partial(*)}(Y)):\beta < \mu_{2,2}\}$
\newline
(this is an overkill).
\end{enumerate}
\end{enumerate}
\end{claim}

\begin{PROOF}{\ref{y.21}}
\mn
\begin{enumerate}
\item[$\oplus_1$]  $\xi_2$ is an ordinal.
\end{enumerate}
\mn
[Why?   To prove that $\xi_2$ is an ordinal we have to assume
$\alpha < \beta \in \xi_2$ and prove $\alpha \in \xi_2$.  As $\beta
\in \xi_2$ clearly $\beta \in \Xi_{\gy}[c \ell]$ for some ${\gy}
\in \text{ Fil}^4_{\partial(*)}(Y,\mu_1)$ for which $D^{\gy}_1 = D$ so there is
$f \in {}^Y(\xi_1)$ such that $f  \rest Z^{\gy} \in {\cF}_{{\gy},\beta}$.  So
rk$_{D+Z[\gy]}(f) = \beta$ hence by \ref{0.A} there is $g
\in {}^Y \lambda$ such that $g \le f$, i.e., $(\forall y \in Y)(g(y)
\le f(y))$ and rk$_{D+Z[\gy]}(g) = \alpha$.  By \ref{r.4}(4) there is 
${\gz} \in \text{ Fil}^4_{\partial(*)}(Y,\mu_1)$ such that 
$D^{\gz}_1 = D + Z[\gy]$ 
and $g \rest Z^{\gz} \in {\cF}_{\gz,\alpha}$ so we are done
proving $\xi_2$ is an ordinal.]

We define the function $c \ell'$ with domain $[\xi_2]^{< \hrtg(Y)}$ as  
follows: 
\mn
\begin{enumerate}
\item[$\oplus_2$]  $c \ell'(u) = \{0\} \cup\{\alpha$: there 
is ${\gy} \in 
\text{ Fil}^4_{\partial(*)}(Y,\mu_1)$ such that $f_{{\gy},\alpha}[c
\ell]$ is well defined \footnote{We could have used
$\{t \in Y:f_{\eta,\alpha}[c \ell](t) \in c \ell(\bold v(u))\} 
\ne \emptyset$ mod $D^{\gy}_2$; also we could have added $u$ to $c
\ell'(u)$ but not necessarily by $\boxplus_2$.}
 and Rang$(f_{{\gy},\alpha}[c \ell]) \subseteq c \ell(\bold v[u])\}$.  
\end{enumerate}
\mn
where
\mn
\begin{enumerate}
\item[$\oplus_3$]   $\bold v[u] := \cup \{c \ell(v):v \subseteq \xi_1$ is of
cardinality $\le \partial$ and is $\subseteq \bold w(v)\}$ .
\end{enumerate}
\mn
where
\mn
\begin{enumerate}
\item[$\oplus_4$]   for $v \subseteq \xi_1$ we let $\bold w(v) = 
\cup\{\text{Rang}(f_{{\gz},\beta}[c \ell]):{\gz} \in 
\text{ Fil}^4_{\partial(*)}(Y,\mu_1)$ and $\beta \in u$ and
$f_{{\frak z},\beta}[c \ell]$ is well defined$\}$.  
\end{enumerate}
\mn
Note that
\mn
\begin{enumerate}
\item[$\oplus_5$]   $c \ell'(u) = \{0\} \cup \{\text{\rm rk}_D(f):
D \in \text{\rm Fil}_{\partial(*)}(Y),Z \in D^+$ 
and $f \in {}^Y \bold v(u)\}$.
\end{enumerate}
\mn
Note that (by \ref{r.4}(1)):
\mn
\begin{enumerate}
\item[$\boxtimes_1$]    for each $u \subseteq \xi_1$ and
${\gx} \in \text{ Fil}^4_{\partial(*)}(Y,\mu_1)$  the set 
$\{\alpha < \xi_2:f_{{\gx}, \alpha}[c \ell]$ is a well defined 
function into $u\}$ has cardinality $< \wlor(T_{D^{\gy}_2}(u))$, that is,
$\langle f_{{\gx},\alpha}[c \ell]:\alpha \in \Xi_{\gx} \cap
\xi_2\rangle$ is a sequence of functions from $Z^{\gx}$ to $u
\subseteq \xi_1$, any two are equal only on a set $= \emptyset$ mod
$D^{\gx}_2$ (with choice it has cardinality $\le {}^{|Y|}|u|$)), call this 
bound $\mu'_{|u,\gx|}$.
\end{enumerate}
\mn
Note
\mn
\begin{enumerate}
\item[$\boxtimes_2$]   if $u_1 \subseteq u_2 \subseteq \xi_2$ \then
\sn
\begin{enumerate}
\item[$(\alpha)$]   $\bold w(u_1) \subseteq \bold w(u_2)$ and
$\bold v(u_1) \subseteq \bold v(u_2) \subseteq \xi_1$
\sn
\item[$(\beta)$]   $c \ell'(u_1) \subseteq c \ell'(u_2)$
\sn
\item[$(\gamma)$]   $u \subseteq \bfv(u)$ and $\bfw[u] \subseteq
  \bfv[u]$
\sn
\item[$(\delta)$]   $u_1 \subseteq c \ell'(u_1)$.
\end{enumerate}
\end{enumerate}
\mn
[Why?  E.g. for clause $(\delta)$; assume $\alpha \in u$ and let $f$
be a unique function from $Y$ into $\{\alpha\}$.  Hence for some $\gy
\in \Fil^4_{\partial(*)}(Y,\mu_1)$ we have $f_{\gy,\alpha}$ is well defined.
Now $\Rang(f_{\gy,\alpha}) \subseteq \bfw(u)$ by the choice of
$\bfw(u)$ in $\oplus_4$ and so $\Rang(f_{\gy,\alpha}) \subseteq \bfv(u)$
  by clause $(\gamma)$ of $\boxplus_2$ hence $\Rang(f_{\gy,\alpha})
  \subseteq c \ell(\bfv,u)$ by the assumption on $c \ell$, see by
  \ref{y.21}(a),(b) and \ref{22.1}(2).  So we have $f_{\gy,\beta}$
 well defined and $\Rang(f_{\gy,\alpha}) \subseteq c \ell(\bfv(u))$ so
    by the definition of $c \ell'(u)$ in $\oplus_2$ we have $\alpha
    \in c \ell'(u)$ so we are done.]
\mn
\begin{enumerate}
\item[$\boxtimes_3$]   if $u \subseteq \xi_2,|u| <$ \hrtg$(Y)$
then $\bold w(u) = \{f_{{\gy},\alpha}(z):\alpha \in u,{\gy}\in \text{
Fil}^4_{\partial(*)}(Y,\mu_1),f_{\gy,\alpha}$ is well defined 
and $z \in Z^{\gy}\}$ is a subset
of $\xi_1$ of cardinality $<$ \hrtg$(|u| \times 
\text{ Fil}^4_{\partial(*)}(Y,\mu_1)) \le \sup\{\hrtg(\beta) \times
\text{ Fil}^4_{\partial(*)}(Y,\mu_1)):\beta < \hrtg(Y)\}$ 
which was named $\mu_{2,1}$ in \ref{y.21}$(e)(\beta)$
\sn
\item[$\boxtimes_4$]   if $u \subseteq \xi_1$ and $|u| < \mu_{2,1}$
then $\cup\{c \ell(v):v \in [u]^{\le \partial}\}$ is a subset of
$\mu_1$ of cardinality $< \hrtg(\mu_1 \times [u]^{\le \partial}) 
\le \sup_{\alpha < \mu_{2,1}} \hrtg(\mu_1 \times [\alpha]^{\le
\partial})$ which we call $\mu_{2,2}$ in \ref{y.21}$(e)(\gamma)$
\sn
\item[$\boxtimes_5$]   if $u \subseteq \xi_2$ and $|u| < \hrtg(Y)$
then $\bold v(u)$ has cardinality $< \mu_{2,2}$.
\end{enumerate}
\mn
[Why?  By $\oplus_3$ and $\boxtimes_3$ and $\boxtimes_4$.]
\mn
\begin{enumerate}
\item[$\boxtimes_6$]   if $u \subseteq \xi_2$ and $|u| < \hrtg(Y)$
then $c \ell'(u) \subseteq \xi_2$ and has cardinality $< \mu_{2,3}$ is
defined in \ref{y.21}$(e)(\delta)$
which we call $\mu_2$.
\end{enumerate}
\mn
[Why?  Without loss of generality 
$\bold v(u) \ne \emptyset$.  By $\oplus_5$ we have
$|c \ell'(u)| < \hrtg({}^Y \bold v(u)) \times \Fil_{\partial(*)}(Y))$
 and by $\boxplus_5$ the latter is $\le
\sup\{\hrtg({}^Y \beta \times \text{ Fil}_{\partial(*)}(Y)):\beta <
\mu_{2,2}\} = \mu_{2,3}$ recalling clause $(e)(\delta)$ of the claim,
so we are done.]
\mn
\begin{enumerate}
\item[$\boxtimes_7$]   $c \ell'$ is a very weak closure operation on
$\lambda$ and has character $(< \mu_2,\hrtg(Y))$.
\end{enumerate}
\mn
[Why?  In Definition \ref{0.F}(1), clause (a) holds by the Definition
of $c \ell'$, clause (b) holds by $\boxplus_6$ and as for clause (c),
$0 \in c \ell'(u)$ by the definition of $c \ell'$ and $u \subseteq c
\ell'(u)$ by clause $(\delta)$ of $\boxtimes_2$.]

Now it is enough to prove
\mn
\begin{enumerate}
\item[$\boxtimes_8$]   $c \ell'$ witnesses Ax$^0_{\xi_2,< \mu_2,Y}$.
\end{enumerate}
\mn
Recalling $\boxtimes_7$, toward contradiction assume $\bar{\cU} = \langle 
{\cU}_\varepsilon:\varepsilon < \partial\rangle$ is 
$\subseteq$-decreasing, ${\cU}_\varepsilon \in [\xi_1]^{< \hrtg(Y)}$  
and $\varepsilon < \partial \Rightarrow {\cU}_\varepsilon 
\nsubseteq c \ell({\cU}_{\varp +1})$.  We define
$\bar \gamma = \langle \gamma_\varepsilon:\varepsilon < \partial\rangle$ by

\[
\gamma_\varepsilon = \text{ Min}({\cU}_\varepsilon \backslash c \ell
({\cU}_{\varepsilon +1})).
\]

\mn
As AC$_\partial$ follows from DC$_\partial$, we can choose $\langle
{\gy}_\varepsilon:\varepsilon < \partial \rangle$ such that
$f_{{\gy}_\varepsilon,\gamma_\varepsilon}[c \ell]$ is well
defined for $\varepsilon < \partial$.

Let for $\varepsilon < \partial$

\[
u_\varepsilon = \{\gamma_\zeta:\zeta \in [\varepsilon,\partial)\}.
\]

\mn
So
\mn
\begin{enumerate}
\item[$(*)_1$]   $u_\varepsilon \in [\xi_1]^{\le \partial} \subseteq
[\xi_1]^{< \hrtg(Y)}$.
\end{enumerate}
\mn
[Why?  By clause (a) of the assumption of \ref{y.21}.]
\mn
\begin{enumerate}
\item[$(*)_2$]   $u_\varepsilon$ is $\subseteq$-decreasing with $\varepsilon$.
\end{enumerate}
\mn
[Why?  By the definition.]
\mn
\begin{enumerate}
\item[$(*)_3$]  $\gamma_\varepsilon \in u_\varepsilon \backslash c
\ell(u_{\varepsilon +1})$ for $\varepsilon < \partial$.
\end{enumerate}
\mn
[Why?  $\gamma_\varepsilon \in u_\varepsilon$ by the definition of
$u_\varepsilon$.]

Now if $\zeta \in [\varp,\gamma)$ then $f_{\gy_\zeta,\gamma_\zeta}[c
\ell]$ is well defined and $\gamma_\zeta \in \cU_\zeta \backslash c
\ell(\cU_{\zeta +1})$ (see the choice of $\gamma_\varp$) but $\langle
\cU_\xi:\xi < \partial\rangle$ is $\subseteq$-decreasing hence
$\gamma_\zeta \in \cU_\zeta$, by the definition of
$\bfw[u_\varp],\Rang(f_{\gy_\zeta,\gamma_\zeta}) \in \bfw(\cU_\varp)$,
hence $\Rang(f_{\gy_\zeta,\gamma_\zeta}) \in \bfv(\cU_\varp) \subseteq
c \ell(\bfv(\cU_\varp))$.  As this holds for every $\zeta \in
[\varp,\gamma)$ we can deduce $u_varp = \{\gamma_\zeta:\zeta \in
[\varp,\partial)\} \subseteq c \ell'(\bfv(\cU_\varp))$.

Lastly, $\gamma_\varepsilon \notin \bold v({\cU}_{\varepsilon +1})$
by the choice of $\beta_\varepsilon$. So $\langle
u_\varepsilon:\varepsilon < \partial \rangle$ contradict the
assumption on $(\xi_1,c \ell)$.
\sn
\relax From the above the conclusion should be clear.  
\end{PROOF}

\begin{claim}
\label{r.4A}
Assume $\kappa < \kappa = 
{ \text{\rm cf\/}}(\lambda) < \lambda$ hence $\kappa$ is regular
$\ge\partial$ of course, and $D$ is the club filter on $\kappa$ 
and $\bar\lambda = \langle \lambda_i:i < \kappa\rangle$ is 
increasing continuous with limit $\lambda$.

Then $\lambda^+ \le \{\text{\rm rk}_{D_\kappa}(f):f \in \prod\limits_{i
< \kappa^+} \lambda^+_i\}$.
\end{claim}

\begin{PROOF}{\ref{r.4A}}
For each $\alpha < \lambda^+$ there is a one to one
\footnote{but, of course, possibly
there is no such sequence $\langle f_\alpha:\alpha < \lambda^+ \rangle$}
 function $g$ from $\alpha$ into $|\alpha| \le \lambda$ and we let
$f \in \prod\limits_{i < \kappa} \lambda_i$ be

\[
f(i) = \text{ otp}(\{\beta < \alpha:g(\beta) < \lambda_i\}.
\]

\mn
Let 

\begin{equation*}
\begin{array}{clcr}
{\cF}_\alpha = \{f: &f \text{ is a function with domain } \kappa
  \text{ satisfying } i < \kappa \Rightarrow f(i) < \lambda^+_i \\
  &\text{ such that for some one to one function } g \text{ from }
  \alpha \text{ into } \lambda \\
  &\text{ for each } i < \kappa \text{ we have } f(i) 
= \text{ otp}(\{\beta < \alpha:g(\beta) < \lambda_i\})\}.
\end{array}
\end{equation*}

\mn
Now
\mn
\begin{enumerate}
\item[$(*)_1$]
\begin{enumerate}
\item[$(\alpha)$]  $\cF_\alpha \ne \emptyset$ for $\alpha < \lambda^+$
\sn
\item[$(\beta)$]  $\langle {\cF}_\alpha:\alpha <
\lambda^+\rangle$ exists as it is well defined
\end{enumerate}
\end{enumerate}
\mn
[Why?  For clause $(\alpha)$ let $g:\alpha \rightarrow \lambda$ be one
to one and so the $f$ defined above belongs to $\cF_\alpha$.  For clause
$(\beta)$ see the definition of $\cF_\alpha$ (for $\alpha <
\lambda^+$).]
\mn
\begin{enumerate}
\item[$(*)_2$]
\begin{enumerate}
\item[$(\alpha)$]   if $f \in {\cF}_\beta, \alpha <
\beta < \lambda^+$ \then \, for some $f' \in \cF_\alpha$ we have $f'
<_{J^{\text{bd}}_\kappa} f$
\sn
\item[$(\beta)$]  $\langle \min\{\text{rk}_D(f):f \in
{\cF}_\alpha\}:\alpha < \lambda^+\rangle$ is strictly increasing hence
$\min\{\rk_D(f):f \in \cF_\alpha\} \ge \alpha$. 
\end{enumerate}
\end{enumerate}
\mn
[Why?  For clause $(\alpha)$, let $g$ witness ``$f \in \cF_\beta$" and
define the function $f' \in \prod\limits_{i < \kappa} \lambda^+_i$ by
$f'(i) = \text{ otp}\{\gamma < \alpha:g(\gamma) < \lambda_i\}$.  So $g
\rest \alpha$ witness $f' \in \cF_\alpha$, and letting $i(*) = \text{
min}\{i:g(\alpha) < \lambda_i\}$ we have $i \in [i(*),\kappa)
\Rightarrow f'(i) < f(i)$ hence $f' <_{J^{\text{bd}}_\kappa} f$ as
promised.  For clause $(\beta)$ it follows.]

Note that
\mn
\begin{enumerate}
\item[$(*)_3$]  if $f \in {\cF}_\alpha$ then, for part (2), for some
${\gy} \in \text{ Fil}^4_{\partial(*)}(\lambda,c \ell)$ and $\beta \ge
\alpha$ we have $f \rest Z[\gy] \in \cF_{\gy,\beta}$.
\end{enumerate}
\mn
[Why?  By $(*)_1 + (*)_2$.]

So we have proved \ref{r.4A}.
\end{PROOF}

\begin{conclusion}
\label{r.5}
1) Assume
\mn
\begin{enumerate}
\item[$(a)$]    {\rm Ax}$^0_{\lambda,< \mu,\kappa}$
\sn
\item[$(b)$]   $\lambda > \text{\rm cf}(\lambda) 
= \kappa$ (not really needed in part (1)).
\end{enumerate}
\mn
\Then \, for some ${\cF}_* \subseteq {}^\kappa \lambda =:
\{f:f$ a partial function from $\kappa$ to $\lambda\}$ we have
\mn
\begin{enumerate}
\item[$(\alpha)$]   every $f \in {}^\kappa \lambda$ is a countable
union of members of ${\cF}_*$
\sn
\item[$(\beta)$]   ${\cF}_*$ is the union of
$|\text{\rm Fil}^4_{\partial(*)}(\kappa,< \mu)|$ well 
ordered sets: $\{{\cF}^*_{\gy}:{\gy} 
\in \text{\rm Fil}^4_{\partial(*)}(\kappa,\mu)\}$
\sn
\item[$(\gamma)$]   moreover there is a function giving for each
$\gy \in \text{\rm Fil}^4_{\partial(*)}(\kappa)$ a well ordering of
${\cF}^*_{\gy}$.
\end{enumerate}
\mn
2) Assume in addition that \hrtg$(\text{\rm Fil}^4_{\partial(*)}(\kappa,<
\mu)) < \lambda,\cf(\lambda^+)$ and \hrtg$({}^\kappa \mu) < \lambda$ 
\then \,  for some ${\gy} \in \text{\rm Fil}^4_{\partial(*)}(\kappa)$
we have $|\cF^*_{\gy}| > \lambda$.

\noindent
3) If in part (2) we omit the assumption on $\cf(\lambda^+)$ still
$\lambda^+ = \sup\{\otp(\Xi_{\gy} \cap \lambda^+):\gy \in
\Fil^4_{\partial()*}(\kappa,\mu)\}$. 
\end{conclusion}

\begin{PROOF}{\ref{r.5}}
1) By the proof of \ref{r.2}.

\noindent
2) Assume that this fails; so for every $\gy \in
\Fil^4_{\partial(*)}(\kappa,< \mu)$, the set $S_{\gy} = \Xi_{\gy}
\cap \lambda^+$ has order type $< \lambda^+$.  But we are assuming
$\cf(\lambda^+) \ge \hrtg(\Fil^4_{\partial()*}(\kappa,\mu)$, so there
is $\gamma < \lambda^+$ such that $\gamma > \otp(S_{\gy})$ for every
relevant $\gy$, \wilog \, $\gamma > \lambda$ and let $g$ be a
one-to-one function from $\gamma$ onto $\lambda$.

We choose $f \in {}^\kappa \lambda$ by

\begin{equation*}
\begin{array}{clcr}
f(i) = \text{ Min}(\lambda \backslash \{f_{{\gy},\alpha}(i): &{\gy} 
\in \text{ Fil}^4_{\partial(*)}(\kappa,\mu) \\
  &f_{\gy,\alpha}(i) \text{ is well defined, i.e.} \\
  &i \in Z[{\gy}] \text{ and } \alpha \in \Xi_{\gy} \text{ and} \\
  &g(\otp(\alpha \cap \Xi_{\gy})) < \mu_i\}).
\end{array}
\end{equation*}

\mn
Now $f(i)$ is well defined as the minimum is taken over a non-empty
set of ordinals, this holds as we substruct from $\lambda$ a set 
which has cardinality $\le \mu_i$ which is $< \lambda$.  
But $f$ contradicts part (1).  Note that in fact $f \in
\prod\limits_{i} \mu^+_i$.

\noindent
3) Same proof as in part (2).
\end{PROOF}

\begin{conclusion}
\label{r.6}
Assume Ax$^0_{\lambda,< \mu,\kappa}$ so $\lambda > \mu$.

\noindent
$\lambda^+$ is not measurable (even in cases it is 
regular\footnote{the regular holds many times by \ref{s.2}}) 
\when \,
\mn
\begin{enumerate}
\item[$\boxtimes$]   $(a) \quad \lambda > \text{\rm cf}(\lambda) = 
\kappa > \aleph_0$
\sn
\item[${{}}$]   $(b) \quad \lambda > 
\hrtg((\text{\rm Fil}^4_{\partial(*)} (\kappa,\mu))$. 
\end{enumerate}
\end{conclusion}

\begin{PROOF}{\ref{r.6}}
Naturally we fix a witness $c \ell$ for Ax$^0_{\lambda,<\mu,\kappa}$.
Let ${\cF}_{\gy},\Xi_{\gy},f_{{\gy},\alpha},{\cF}^\lambda_{{\gy},\alpha}$ be 
defined as in \ref{r.4} so by claims \ref{r.4}, \ref{r.4A} we have 
$\cup\{\Xi_{\gy}:{\gy} \in \text{ Fil}^4_{\partial(*)}(\kappa)\} 
\supseteq \lambda^+$; moreover, $\alpha \in \lambda^+ \cap \Xi_{\gy}
\Rightarrow f_{\eta,\alpha} \in {}^\kappa \lambda$.

Let ${\gy} \in \text{ Fil}^4_{\partial(*)}(\kappa,\mu)$ be
such that $|{\cF}_{\gy}| > \lambda$, we can find such $\gy$ by \ref{r.5},
 as \wilog \, we can assume $\lambda^+$ is regular (or even
measurable, toward contradiction).  Let $Z = Z[{\gy}]$.  So $\Xi_{\gy}$ is 
a set of ordinals of cardinality $>\lambda$.  For $\zeta <\text{ otp}
(\Xi_{\gy})$ let $\alpha_\zeta$ be the $\zeta$-th member of 
$\Xi_{\gy}$, so $f_{{\gy},\alpha_\zeta}$ is well defined.  Toward contradiction
let $D$ be a (non-principal) ultrafilter on $\lambda^+$ which is 
$\lambda^+$-complete.  For $i \in Z$ let $\gamma_i <
\lambda$ be the unique ordinal $\gamma$
such that $\{\zeta < \lambda^+:f_{{\gy},\alpha_\zeta}(i)
= \gamma\} \in D$.  As $|Z| \le \kappa < \lambda^+$ and $D$ is
$\kappa^+$-complete  clearly $\{\zeta:\bigwedge\limits_{i \in Z}
f_{{\gy},\alpha_\zeta}(i) = \gamma_i\} \in D$, so as $D$ is a 
non-principal ultrafilter, for some $\zeta_1 < \zeta_2,f_{{\gy},
\alpha_{\zeta_1}} = f_{{\gy},\alpha_{\zeta_2}}$, contradiction.   
So there is no such $D$.
\end{PROOF}

\begin{remark}  
\label{r.7n}
Similarly if $D$ is $\kappa^+$-complete and
weakly $\lambda^+$-saturated and Ax$^0_{\lambda^+,< \mu}$, see \cite{Sh:1005}.
\end{remark}

\begin{claim}
\label{r.9}
If {\rm Ax}$^0_{\lambda,< \mu,\kappa}$, 
\then \, we can find $\bar C$ such that:
\mn
\begin{enumerate}
\item[$(a)$]   $\bar C = \langle C_\delta:\delta \in S \rangle$
\sn
\item[$(b)$]   $S = \{\delta < \lambda:\delta$ is a limit ordinal of
cofinality $\ge \partial(*)\}$
\sn
\item[$(c)$]   $C_\delta$ is an unbounded subset of $\delta$, even a club
\sn
\item[$(d)$]   if $\delta \in S$, {\rm cf}$(\delta) \le \kappa$ then
$|C_\delta| < \mu$
\sn
\item[$(e)$]   if $\delta \in S$, {\rm cf}$(\delta) > \kappa$ then
$|C_\delta| < \hrtg(\mu \times [\text{\rm cf}(\delta)]^\kappa)$.
\end{enumerate}
\end{claim}

\begin{remark}  
1) Recall that if we have Ax$^4_\lambda$ (see \ref{22.1}(5)) then
trivially there is $\langle C_\delta:\delta < \lambda$, cf$(\delta)
\le \partial\rangle,C_\delta$ a club of $\delta$ of order type
cf$(\delta)$ as if $<_*$ well order $[\lambda]^{\le\partial}$ we let
$C_\delta :=$ be the $<_*$-minimal $C$ which is a closed unbounded
subset of $\delta$ of order type cf$(\delta)$.

\noindent
2) Ax$^0_{\lambda,< \xi,\kappa}$ suffices 
if $\kappa < \xi < \lambda$.
\end{remark}

\begin{PROOF}{\ref{r.9}}
The ``even a club" is not serious as we can replace
$C_\delta$ by its closure in $\delta$.

Let $c \ell$ witness Ax$^0_{\lambda,< \mu,\kappa}$.  For each $\delta
\in S$ with cf$(\delta) \in [\partial(*),\kappa]$ we let

\[
C_\delta = \cap\{\delta \cap c \ell(C):C \text{ a club of } \delta
\text{ of order type cf}(\delta)\}.
\]

\mn
Now $\bar C'= \langle C_\delta:\delta \in S$ and cf$(\delta) 
\in [\partial(*),\kappa] \rangle$ is well defined and exist.  
Clearly $C_\delta$ is a subset of $\delta$.

For any club $C$ of $\delta$ of order type cf$(\delta) \in
[\partial(*),\kappa]$ clearly $\delta \cap c \ell(C) \subseteq c \ell(C)$
which has cardinality $< \mu$. 

The main point is to show that $C_\delta$ is unbounded in $\delta$, otherwise
we can choose by induction on $\varepsilon <
\partial$, a club $C_{\delta,\varepsilon}$ of $\delta$ of order type
cf$(\delta)$, decreasing with $\varepsilon$ such that
$C_{\delta,\varepsilon} \nsubseteq c \ell(C_{\delta,\varepsilon +1})$,
we use DC$_\partial$.  But this contradicts the choice of $c \ell$
recalling Definition \ref{22.1}(1).

If $\delta < \lambda$ and cf$(\delta) > \kappa$ we let

\begin{equation*}
\begin{array}{clcr}
C^*_\delta = \cap\{\cup\{\delta \cap c \ell(u): &u \subseteq C \text{ has
cardinality } \le \partial\}: \\
  &C \text{ is a club of } \delta \text{ of order type cf}(\delta)\}.
\end{array}
\end{equation*}

\mn
A problem is a bound of $|C^*_\delta|$.  Clearly for $C$ a
club of $\delta$ of order type cf$(\delta)$ the order-type of the set
$\cup \{\delta \cap c \ell(v):v
\subseteq C$ has cardinality $\le \partial\}$ is $< \hrtg(\mu \times
[\text{cf}(\delta)]^\kappa)$.  As for ``$C^*_\delta$ is a club" it is
proved as above.
\end{PROOF}

\noindent
The following lemma gives the existence of a class of regular
successor cardinals.
\begin{lemma}
\label{s.2}
1) Assume
\mn
\begin{enumerate}
\item[$(a)$]   $\delta$ is a limit ordinal $< \lambda_*$ with {\rm
cf}$(\delta) = \partial$
\sn
\item[$(b)$]    $\lambda^*_i$ is a cardinal for $i < \delta$
increasing with $i$
\sn
\item[$(c)$]   $\lambda_* = \Sigma\{\lambda^*_i:i < \delta\}$
\sn
\item[$(d)$]   $\lambda^*_{i+1} \ge \hrtg(\mu \times
{}^\kappa(\lambda^*_i))$ for $i < \delta$ and
$(\alpha) \vee (\beta)$ hold where:
\sn
\begin{enumerate}
\item[$(\alpha)$]  Ax$^4_\lambda$ \underline{or}
\sn
\item[$(\beta)$]   $\lambda^*_{i+1} \ge \hrtg
({\text{\rm Fil\/}}^4_{\partial(*)}(\lambda^*_i,\mu))$ and
\hrtg$([\lambda^*_i]^{\le \kappa}) \le \lambda^*_{i+1}$
\end{enumerate}
\item[$(e)$]   {\rm Ax}$^0_{\lambda,<\mu,\kappa}$ and $\mu < \lambda^*_0$
\sn
\item[$(f)$]   $\lambda = \lambda^+_*$
\end{enumerate}
\mn
\Then \, $\lambda$ is a regular cardinal.

\noindent
2) Assume {\rm Ax}$^4_\lambda,\lambda = \lambda^+_*,\lambda_*$
singular and $\chi < \lambda_* \Rightarrow \hrtg({}^\partial \chi) \le
\lambda_*$ \then \, $\lambda$ is regular.
\end{lemma}

\begin{remark}  
This says that the successor of many strong limit singulars is regular.
\end{remark}

\begin{question}
1) Is \hrtg$({\cP}({\cP}(\lambda^*_i))) \ge$ 
\hrtg$(\text{Fil}^4_{\partial(*)}(\lambda^*_i))$?

\noindent
2) Is $|c \ell(f \restriction B)| \le$ \hrtg$([B]^{< \aleph_0})$ for the
natural $c \ell$ and $f,B$ as in the proof of \ref{s.2}?
\end{question}

\begin{PROOF}{\ref{s.2}}
1) We can replace $\delta$ by cf$(\delta)$ so \wilog \,
$\delta$ is a regular cardinal so $\delta = \partial$.

So 
\mn
\begin{enumerate}
\item[$(*)_1$]  
\begin{enumerate}
\item[(a)]  fix $c \ell:[\lambda]^{\le \kappa} \rightarrow {\cP}(\lambda)$
a witness to Ax$^0_{\lambda,< \mu,\kappa}$
\sn
\item[(b)]   let $\langle C_\xi[c \ell]:\xi < \lambda$, cf$(\xi) \ge 
\partial \rangle$ be as in the proof of \ref{r.9}, so $\xi < \lambda
\wedge \partial \le \text{ cf}(\xi) < \lambda \Rightarrow |C_\xi[c
\ell]| <  \lambda$.
\end{enumerate}
\end{enumerate}
\mn
[Why the last inequality?  If $\delta < \lambda^+$, then there is $i$ such
that $\lambda^*_i > \mu + \cf(\partial)$ hence $\otp(C_\delta) <
\hrtg(\mu \times [\cf(\delta)]^\kappa) \le \hrtg([\lambda^*_i]^\kappa)
< \lambda^*_{i+1}$.]

First, we shall use just $\lambda > \lambda_* \wedge (\forall \delta <
\lambda)(\text{\rm cf}(\delta) < \lambda_*)$, a weakening of the 
assumption that $\lambda = \lambda^+_*$.  

Now
\mn
\begin{enumerate}
\item[$\boxtimes_1$]    for every $i < \delta$ and
$A \subseteq \lambda$ of cardinality $\le \lambda^*_i$, we can find
$B \subseteq \lambda$ of cardinality $\le \lambda_*$ satisfying
$(\forall \alpha \in A)$[$\alpha$ is limit $\wedge \text{
cf}(\alpha) \le \lambda^*_i \Rightarrow \alpha = \sup(\alpha \cap B)]$.
\end{enumerate}
\mn
The proof of this will take some time.  By \ref{r.9} (and \ref{x.5}) 
the only problem is for 
$Y := \{\alpha:\alpha \in A,\alpha > \sup(A \cap \alpha),\alpha$ 
a limit ordinal of cofinality $< \partial +
\aleph_1\}$; so $|Y| \le \lambda^*_i$.  Note: if we assume
Ax$^4_\lambda$ this would be immediate.

We define $D$ as the family of sets $A \subseteq Y$ such that:
\mn
\begin{enumerate}
\item[$\circledast^1_A$]   for some set $C \subseteq \lambda$ of $\le
\partial$ ordinals, the set $B_C =: \cup\{\text{Rang}(f_{{\gx},
\zeta}):{\gx} \in \text{Fil}^4_{\partial(*)}(\lambda^*_i,\mu)$ and 
$\zeta \in C$ \underline{or} for some $\xi \in C$, we have
$\lambda^*_i \ge \text{ cf}(\xi) > \partial$ 
and $\zeta \in C_\xi[c \ell]\}$ satisfies 
$\alpha \in Y \backslash A \Rightarrow \alpha = \sup(\alpha \cap B_C)$. 
\end{enumerate}
\mn
Clearly
\mn
\begin{enumerate}
\item[$\circledast_2$]   $(a) \quad Y \in D$
\sn
\item[${{}}$]   $(b) \quad D$ is upward closed
\sn
\item[${{}}$]   $(c) \quad D$ is closed under intersection of $\le
\partial$ hence of $< \partial(*)$ sets.
\end{enumerate}
\mn
[Why?  For clause (a) use $C = \emptyset$, for clause (b), note that if $C$
witness a set $A \subseteq Y$ belongs to $D$ then 
it is a witness for any $A' \subseteq Y$ such that $A \subseteq A'$.  
Lastly, for clause (c) if $A_\varepsilon \in D$ for $\varepsilon <
\varepsilon(*) <\partial^+$, as we have AC$_\partial$, there is a
sequence $\langle C_\varepsilon:\varepsilon < \varepsilon(*)\rangle$
such that $C_\varepsilon$ witnesses $A_\varepsilon \in D$ for
$\varepsilon < \varepsilon(*) < \partial^+$, then
$C := \cup\{C_\varepsilon:\varepsilon < \varepsilon(*)\}$ witnesses 
$A := \cap \{A_\varepsilon:\varepsilon < \varepsilon(*)\} \in D$ and,
again by AC$_\partial$, we have $|C| \le \partial$.]
\mn
\begin{enumerate}
\item[$\circledast_3$]   if $\emptyset \in D$ then we are done.
\end{enumerate}
\mn
[Why?  For $a = \emptyset \in D$ let $C \subseteq \lambda$ be as
promised in $\circledast_1$ and then $B_C$ is as required; 
its cardinality $\le \lambda^*_{i+1}$ by \ref{r.9}.]

So assume $\emptyset \notin D$, so $D$ is an $\partial^+$-complete
filter on $Y$.  As $1 \le |Y| \le \lambda^*_i$, let $g$ be a one to
one function from $|Y| \le \lambda^*_i$ onto $Y$ and let
\mn
\begin{enumerate}
\item[$\circledast_4$]  
\begin{enumerate}
\item[(a)]  $D_1 := \{B \subseteq \lambda^*_i:
\{g(\alpha):\alpha \in B \cap |Y|\} \in D\}$
\sn
\item[(b)]  $\zeta := \text{\rm rk}_{D_1}(g)$
\sn
\item[(c)]  $D_2 := \{B \subseteq \lambda^*_i:B \notin D_1$ and 
$\text{\rm rk}_{D_1+(\lambda^*_i \backslash B)}(g) > \zeta\} \cup D_1$. 
\end{enumerate}
\end{enumerate}
\mn
So $D_2$ is an $\partial^+$-complete filter on $\lambda^*_i$ extending
$D_1$.

Let $B_* \in D_2$ be such that $(\forall B')[B' \in D_2 \wedge B'
\subseteq B_* \Rightarrow c \ell(\text{Rang}(g \restriction B')) \supseteq
(\text{Rang}(g \restriction B_*)]$.  Let ${\cU} = \cap\{c
\ell(\text{Rang}(g \restriction B'):B' \in D_2\}$, so
Rang$(g \restriction B_*) \subseteq {\cU}$, even equal.

Let $h$ be the function with domain $B_*$ 
defined by $\alpha \in B_* \Rightarrow
h(\alpha) = \text{ otp}(g(\alpha) \cap {\cU})$.

So ${\gx} := (D_1,D_2,B_*,h) \in \text{ Fil}^4_{\partial(*)}
(\lambda^*_i,\mu)$ and for some $\zeta$ we have $g \restriction B_* 
= f_{{\gx},\zeta}[c \ell]$.

It suffices to consider the following two subcases.
\bigskip

\noindent
\underline{Subcase 1a}: cf$(\zeta) > \partial$.

So recalling $(*)_1(b), \, C_\zeta[c \ell]$ is 
well defined and let $C := \{\zeta\}$ hence $B_C =
\cup\{\text{Rang}(f_{\gx,\varepsilon}[c \ell]:\varepsilon \in
C_\zeta[c \ell]\}$ so $C$ exemplifies that the set $X := \{\alpha \in
Y:\alpha > \sup(\alpha \cap B_C)\}$ belongs to $D$ hence $X_* =
\{\alpha < |Y|:g(\alpha) \in X\}$ belongs to $D_1$.

Now define $g'$, a function from $\lambda^*_i$ to Ord by $g'(\alpha) =
\sup(g(\alpha) \cap B_C)+1$ if $\alpha \in X_*$ and $g'(\alpha) = 0$
otherwise.  Clearly $g' < g$ mod $D_1$ hence {\rm rk}$_{D_1}(g') < \zeta$,
hence there is $g'',g' <_{D_1} g'' <_{D_1} g$ such that $\xi :=
\text{\rm rk}_{D_1}(g'') \in C_\zeta[c \ell]$.

Now for some ${\gy} \in \text{ Fil}^4_{\partial(*)}(\lambda^*_i)$ we
have $D^{\gy} = D_2$ and $g'' = f_{{\gy},\xi}$ mod $D^{\gy}_2$.

So $B =: \{\varepsilon < |Y|:g''(\varepsilon) = f_{{\gy},\xi}(\varp)\} \in
D^{\gy}_2$ hence $B \in D^+_2$.  
So $B \cap B_* \cap X_* \in D^+_2$ but if $\varepsilon \in B \cap B_*
\cap A_*$ then $f_{{\gy},\xi}(\varepsilon) \in B_C$ and 
$f_{{\gy},\xi}(\varepsilon) \in \sup((B_C \cap g(\varepsilon)),
g(\varepsilon))$.

This gives contradiction.
\bigskip

\noindent
\underline{Subcase 1b}:  cf$(\zeta) \le \partial$.

We choose a $C \subseteq \zeta$ of order type $\le \partial$ 
unbounded in $\zeta$ and proceed as in subcase 1a.

As we have covered both subcases, we have proved $\boxtimes_1$.

\noindent
Recall we are assuming $\delta = \partial$; now:
\mn
\begin{enumerate}
\item[$\boxtimes_2$]    for every $A \subseteq \lambda$
of cardinality $\le \lambda_*$ there is $B \subseteq \lambda$ of
cardinality $\le \lambda_*$ such that:
\sn
\begin{enumerate}
\item[$\oplus$]  $A \subseteq B,[\alpha +1 \in A \Rightarrow \alpha \in
B]$ and $[\alpha \in A \wedge \aleph_0 \le \text{ cf}(\alpha) <
\lambda_* \Rightarrow \alpha = \sup(B \cap \alpha)]$.
\end{enumerate}
\end{enumerate}
\mn
[Why?  Choose a $\subseteq$-increasing sequence 
$\langle A_j:j < \delta\rangle$ such that $A =
\cup\{A_i:i <\delta\}$ and $j < \delta \Rightarrow
|A_j| \le \lambda^*_j$, possible as $|A|
\le \lambda_*$.  For each $j < \delta$ there exists $B_j$ such that the
conclusion of $\boxplus_1$ holds with $(A_j,B_j,\lambda^*_j)$ here
standing for $(A,B,\lambda_i)$ there, so $|B_j| \le \lambda_*$.  So
as AC$_\delta$ holds (as $\delta
\le \partial$) there is a sequence $\langle \bar B_j:j <
\delta\rangle$, each $\bar B_j$ as above.

Lastly, let $B = \cup\{B_j:j < \delta\}$, it is as required.]
\mn
\begin{enumerate}
\item[$\boxtimes_3$]  for every $A \subseteq \lambda$ of cardinality
$\le \lambda_*$ we can find $B \subseteq \lambda$ of cardinality
$\le \lambda_*$ such that $A \subseteq B,[\alpha +1 \in B
\Rightarrow \alpha \in B]$ and $[\alpha \in B$ is a limit ordinal
$\wedge \text{ cf}(\alpha) < \lambda_* \Rightarrow \alpha = 
\sup(B\cap \alpha)]$.
\end{enumerate}
\mn
[Why?  We choose $B_i$ by induction on $i < \omega \le \partial$ such
that $|B_i| \le \lambda_*$ by $B_0 = A,B_{2i+1} = \{\alpha: 
\alpha \in B_{2i}$ or $\alpha +1 \in
B_{2i+1}\}$ and $B_{2i+2}$ is chosen as $B$ was chosen in
$\boxtimes_2$ for $i$ with $B_{2i+1},B_{2i+2}$ here in the role of
$A,B$ there.  There is such $\langle B_i:i < \omega \rangle$ as $\DC =
\DC_{\aleph_0}$ holds.
So easily $B = \cup\{B_i:i < \omega\}$ is as required.]

Now return to our main case $\lambda = \lambda^+_*$
\mn
\begin{enumerate}
\item[$\boxtimes_4$]    $\lambda^+_*$ is regular.
\end{enumerate}
\mn
[Why?  Otherwise cf$(\lambda^+_*) < \lambda^+_*$ hence
cf$(\lambda^+_*) \le \lambda_*$, but $\lambda_*$ is singular so
cf$(\lambda^+_*) < \lambda_*$ hence
there is a set $A$ of cardinality cf$(\lambda^+_*) < \lambda_*$
such that $A \subseteq \lambda^+_* = \sup(A)$.  Now choose $B$ as in
$\boxtimes_3$.  So $|B| \le \lambda_*,B$ is an unbounded subset of
$\lambda^+_*,\alpha +1 \in B \Rightarrow \alpha \in B$ and if $\alpha
\in B$ is a limit ordinal then cf$(\alpha) \le |\alpha| \le
\lambda_*$, but cf$(\alpha)$ is regular so cf$(\alpha) < \lambda_*$
hence $\alpha =\sup(B \cap \alpha)$.  But this trivially implies that $B =
\lambda^+_*$, but $|B| \le \lambda_*$, contradiction.]

\noindent
2) Similar, just easier.  
\end{PROOF}

\begin{remark}
\label{sp.3}
Of course, if we assume Ax$^4_\lambda$ \then \, the proof of 
\ref{s.2} is much simpler: if $<_*$ is a well
ordering of $[\lambda]^{\le \partial}$ for $\delta < \lambda$ of
cofinality $\le \partial$ let $C_\delta =$ the $<_*$-first closed
unbounded subset of $\delta$ of order type cf$(\delta)$, see \ref{c11}.
\end{remark}

\begin{claim}
\label{sp.4}
Assume
\mn
\begin{enumerate}
\item[$(a)$]   $\langle \lambda_i:i < \kappa \rangle$ is an
increasing continuous sequence of cardinals $> \kappa$
\sn
\item[$(b)$]   $\lambda = \lambda_\kappa = \Sigma\{\lambda_i:i < \kappa\}$
\sn
\item[$(c)$]   $\kappa = \text{\rm cf}(\kappa) > \partial$
\sn
\item[$(d)$]   {\rm Ax}$^0_{\lambda,< \mu,\kappa}$
\sn
\item[$(e)$]  \hrtg$(\text{\rm Fil}^4_{\partial(*)}(\kappa,\mu)) <
\lambda$ and $\kappa,\mu < \lambda_0$ 
\sn
\item[$(f)$]   $S := \{i < \kappa:\lambda^+_i$ is a regular
cardinal$\}$ is a stationary subset of $\kappa$
\sn
\item[$(g)$]  let $D := D_\kappa+S$ 
where $D_\kappa$ is the club filter on $\kappa$
\sn
\item[$(h)$]  $\gamma(*) = \text{\rm rk}_D(\langle \lambda^+_i:i <
\kappa\rangle)$. 
\end{enumerate}
\mn
\Then \, $\gamma(*)$ has cofinality $> \lambda$, so
$(\lambda,\gamma(*)] \cap \text{\rm Reg } \ne \emptyset$.
\end{claim}

\begin{PROOF}{\ref{sp.4}}
Recall \ref{r.4} which we shall use.   Toward
contradiction assume that 
cf$(\gamma(*)) \le \lambda_\kappa$, but $\lambda_\kappa$
is singular hence for some $i(*) < \kappa$, cf$(\gamma(*)) \le
\lambda_{i(*)}$.  Let $c \ell$ witness Ax$^0_{\lambda,< \mu,\kappa}$.

Let $B$ be an unbounded subset of $\gamma(*)$ of order type
cf$(\gamma(*)) \le \lambda_{i(*)}$.  By renaming \wilog \, $i(*)=0$.

For $\alpha < \gamma(*)$ let

\begin{equation*}
\begin{array}{clcr}
{\cU}_\alpha = \cup\{\text{Rang}(f_{{\gy},\alpha}): &f_{{\gy},
\alpha}[c \ell] \text{ is well defined } \in \Pi\{\lambda^+_i:i \in
Z^{\gy}\} \\
  &\text{ and } \gy \in \text{ Fil}^4_{\partial(*)}(\kappa) \text{ and }
  D^{\frak y}_1 = D\}.
\end{array}
\end{equation*}

\mn
Clearly $\cU_\alpha$ is well defined by \ref{r.4}; moreover, $\langle
\cU_\alpha:\alpha < \gamma(*)\rangle$ exists and
$|{\cU}_\alpha| \le$ \hrtg$(\kappa \times 
\text{ Fil}^4_{\partial(*)}(\kappa,\mu)) =$ 
\hrtg$(\text{Fil}^4_{\partial(*)}(\kappa,\mu))$, even $<$ recalling
\ref{x.5}(4).  
Let ${\cU} = \cup\{{\cU}_\alpha:\alpha \in B\}$ so $|{\cU}|
\le \hrtg(\text{Fil}^4_{\partial(*)}(\kappa,\mu)) + |B|$.

We define $f \in \prod\limits_{i < \kappa} \lambda^+_i$ by
\mn
\begin{enumerate}
\item[$(\alpha)$]   $f(i)$ is: $\sup({\cU} \cap \lambda^+_i) +1$ if
cf$(\lambda^+_i) > |{\cU}|$ and zero otherwise.
\end{enumerate}
\mn
So
\mn
\begin{enumerate}
\item[$(\beta)$]   $f \in \prod\limits_{i < \kappa} \, \lambda^+_i$.
\end{enumerate}
\mn
Clearly
\mn
\begin{enumerate}
\item[$(\gamma)$]   $\{i < \kappa:f(i) = 0\} = \emptyset$ mod $D$.
\end{enumerate}
\mn
Let $\alpha(*) = \text{ rk}_D(f)$, it is $< \text{rk}_D(\langle
\lambda^+_i:i < \kappa\rangle) = \gamma(*)$, so by clause $(\gamma)$
there is $\beta(*) \in B$ such that $\alpha(*) < \beta(*) < \gamma(*)$
hence for some $g \in \prod\limits_{i < \kappa} \lambda^+_i$ we have
rk$_D(g) = \beta(*)$ and $f < g$ mod $D$, so for some ${\gy} \in 
\text{ Fil}^4_{\partial(*)}(\kappa)$ we have $D^{\frak y}_1 = D_\kappa + S$ 
and $g \in {\cF}_{\gy,\beta(*)}$, hence $f(i) < g(i) < f_{\gy,\beta(*)}(i) 
\in {\cU} \cap \lambda^+_i$ for every $i \in Z^{\gy} \cap S$.

So we get easy contradiction to the choice of $g$.  
\end{PROOF}

\begin{claim}
\label{sp.6}
Assume $c \ell$ witness {\rm Ax}$^0_{\alpha,< \mu,\kappa}$ and
$\hrtg(Y) < \mu \in [\kappa,\mu)$.  
The ordinals $\gamma_\ell,\ell = 0,1,2$ are nearly
equal see, i.e. $\circledast$ below holds where:
\mn
\begin{enumerate}
\item[$\boxtimes$]
\begin{enumerate}
\item[(a)]  $\gamma_0 = \hrtg({}^Y \alpha)$, a cardinal
\sn
\item[(b)]  $\gamma_1= \cup\{\rk_D(\gamma):\gamma = \rk_D(\alpha)$ for
  some $D \in \text{\rm Fil}_{\partial(*)}(Y)\}$
\sn
\item[(c)]  $\gamma_2 = 
\sup\{\text{\rm otp}(\Xi_{\frak y} [c \ell])+1:{\gy} \in 
\text{\rm Fil}^4_{\partial(*)}(Y)\}$
\end{enumerate}
\sn
\item[$\circledast$]
\begin{enumerate}
\item[$(\alpha)$]  $\gamma_2 \le \gamma_1 \le \gamma_0$
\sn 
\item[$(\beta)$]  $\gamma_0$ is the union of 
{\rm Fil}$^4_{\partial(*)}(Y)$ sets each of order type $< \gamma_2$
\sn 
\item[$(\gamma)$]  $\gamma_0$ is the disjoint union of
$< \hrtg(\cP(\Fil^4_{\partial(*)}(Y)))$ sets each of order type $< \gamma_2$
\sn 
\item[$(\delta)$]  if $\gamma_0 > \hrtg({\cP}
(\text{\rm Fil}^4_{\partial(*)}(Y)))$ and $\gamma_0 \ge
|\gamma_2|^+$ then $|\gamma_0| \le |\gamma_2|^{++}$ 
and {\rm cf}$(|\gamma_2|^+) < \hrtg({\cP}(\text{\rm Fil}^4_{\partial(*)}(Y)))$.
\end{enumerate}
\end{enumerate}
\end{claim}

\begin{PROOF}{\ref{sp.6}}
Straightforward, see \ref{x.5}.
\end{PROOF}
\newpage

\section {Concluding Remarks} \label{3}

In May 2010, David Aspero asked whether it is true that I have results
along the following lines (or that it follows from such a result):

If GCH holds and $\lambda$ is a singular cardinal of uncountable
cofinality, then there is a well-order of $\cH(\lambda^+)$ definable
in $(\cH(\lambda^+),\in)$ using a parameter.

The answer is yes by \cite[4.6,pg.117]{Sh:497} but
we elaborate this below somewhat more generally.  Much earlier Gitik
\cite{Gi80} had proved (using suitable large cardinals) the
consistency of ``$\ZF$ + every infinite cardinal has cofinality
$\aleph_0$, i.e. $\aleph_0$ is the only regular cardinal".  This
naturally raises the question what suffices to have a class of
regulars.  Gitik told me that in Luming 2008 Woodin has conjectured:
\mn
\begin{enumerate}
\item[$\boxplus$]  let $V$ be a model of $\ZF + \DC$, suppose that
  $\kappa$ is a singular strong limit cardinal of cofinality
  $\omega_1$ and $|\cH(\kappa)| = \kappa$.  Is then $\cP(\kappa)$ well
  orderable? 
\end{enumerate}
\mn
Now \cite{Sh:497} gives some information.  The results here (\ref{c2})
confirm $\boxplus$.

\begin{claim}
\label{c2}
{\rm [DC]}
Assume that $\mu$ is a singular cardinal of cofinality $\kappa >
\aleph_0$ (no {\rm GCH} needed), the parameter $X \subseteq \mu$ codes
in particular the tree $\cT = {}^{\kappa >}\lambda$ and the set
$\cP(\cP(\kappa))$ and $F:{}^\omega \mu \rightarrow \mu$ which satisfies
``$(\mu,F)$ has no infinite decreasing $\omega$-chain of subalgebras"; 
in particular, from $X$ a well orderings of $[\lambda]^{< \kappa} \cup
\cP(\cP(\kappa))$ are definable.  \Then \, (with this parameter) 
we can define a well ordering of the set of $\kappa$-branches of the
tree $({}^{\kappa >}\lambda,\triangleleft)$.
\end{claim}

\begin{PROOF}{\ref{c2}}
\underline{Proof of \ref{c2}}:

Let $\langle \text{ cd}_i:i < \kappa\rangle$ satisfies
\mn
\begin{enumerate}
\item[$\boxplus_1$]  cd$_i$ is a one-to-one function from ${}^i \mu$ into
$\mu$, (definable from $X$ uniformly (in $i$))
\sn
\item[$\boxplus_2$]  let $<_\kappa$ be a well ordering of
Fil$^4_\kappa(\kappa)$ definable from $X$.
\end{enumerate}
\mn
For $\eta \in {}^\kappa \mu$ let $f_\eta:\kappa \rightarrow \mu$ be
defined by $f_\eta(i) = \text{ cd}_i(\eta \rest i)$, so $\bar f =
\langle f_\eta:\eta \in{}^\kappa \mu \rangle$ is well defined.

Let $\bar{\cF} = \langle \cF_{\gy}:\gy \in \text{
Fil}^4_\kappa(\kappa)\rangle$ be as in Theorem \ref{r.2} with
$\mu,\kappa$ here standing for $\lambda,Y$ there; there is such
$\bar{\cF}$ definable from $X$ as $X$ codes also $[\mu]^{\aleph_0}$,
see \S1.

So for every $\eta \in {}^\kappa \mu$ there is $\gy \in 
\text{ Fil}^4_\kappa(\kappa)$ such that $f \rest Z_{\gy} \in
\cF_{\gy}$ and $D^{\gy}_1$ contains all co-bounded subsets of $\kappa$ so let
$\gy(\eta)$ be the $<_\kappa$-first such $\gy$.  Now we define a well
ordering $<_*$ of ${}^\kappa \mu$: for $\eta,\nu \in {}^\kappa \mu$
let $\eta <_* \nu$ \Iff \, rk$_{D_1[\gy(\eta)]}(f_\eta \rest
Z_{\gy(\eta)}) < \text{ rk}_{D_1(\gy(\nu))}(f_\nu \rest Z_{\gy(\nu)})$
or equality holds and $\gy(\eta) < \gy(\nu)$.

This is O.K. because
\mn
\begin{enumerate}
\item[$(*)$]  if $\eta \ne \nu \in {}^\kappa \mu$ then $f_\eta(i) \ne
f_\nu(i)$ for every large enough $i < \kappa$ (i.e. $i \ge
\min\{j:\eta(j) \ne \nu(j)\}$.
\end{enumerate}
\end{PROOF}

\begin{conclusion}
\label{c4}
[DC]
Assume $\mu$ is a singular cardinal of uncountable cofinality and
$\cH(\mu)$ is well orderable of cardinality $\mu$ and $X \subseteq
\mu$ codes $\cH(\mu)$ and a well ordering of $\cH(\mu)$.  \Then \, we can 
(with this $X$ as parameter) define a well ordering of $\cP(\mu)$;
hence of $\cH(\mu^+)$.
\end{conclusion}

\begin{PROOF}{\ref{c4}}
\underline{Proof of \ref{c4}}:

Let $\langle \mu_i:i < \kappa\rangle$ be an increasing sequence of
cardinals $< \mu$ with limit $\mu$.  Clearly $2^{\mu_i} < \mu$.

Let $\langle \text{cd}^*_i:i < \kappa\rangle$ satisfies
\mn
\begin{enumerate}
\item[$\boxplus_2$]  $\cd^*_i$ is a one-to-one function from $\cP(\mu_i)$
into $\mu$, (definable uniformly from $X$).
\end{enumerate}
\mn
So cd$_*:\cP(\mu) \rightarrow {}^\kappa \mu$ defined by
$(\text{cd}_*(A))(i) = \text{ cd}^*_i(A \cap \mu_i)$ for $A \subseteq
\mu,i < \kappa$, is a one-to-one function from $\cP(\mu)$ into
${}^\kappa \mu$.  Now use \ref{c2}.
\end{PROOF}

\noindent
We return to \ref{s.2}(2)
\begin{claim}
\label{c11}
{\rm [DC]}
1) The cardinal $\lambda^+$ is regular \when:
\mn
\begin{enumerate}
\item[$\boxplus$]  $(a) \quad$ {\rm Ax}$^4_{\lambda^+}$,
i.e. $[\lambda^+]^{\aleph_0}$ is well orderable
\sn
\item[${{}}$]  $(b) \quad |\alpha|^{\aleph_0} < \lambda$ for $\alpha <
\lambda$
\sn
\item[${{}}$]  $(c) \quad \lambda$ is singular.
\end{enumerate}
\mn
2) Also there is $\bar e = \langle e_\delta:\delta < 
\lambda^+ \rangle,e_\delta \subseteq \delta = \sup(e_\delta),|e_\delta|
\le \text{\rm cf}(\delta)^{\aleph_0}$.
\end{claim}

\begin{remark}
\label{c14}
Compare with \ref{s.2}; we use here more choice, but cover more cardinals.
\end{remark}

\begin{PROOF}{\ref{c11}}
Let $<_*$ be a well ordering of the set $[\lambda^+]^{\aleph_0}$.

As earlier let 
$F:{}^\omega(\lambda^+) \rightarrow \lambda^+$ be such that there
is no $\subset$-decreasing sequence $\langle c \ell_F(u_n):n <
\omega\rangle$ with $u_n \subseteq \lambda^+$.  Let $\Omega = \{\delta
\le \lambda^+:\delta$ a limit ordinal, $\delta < \lambda^+ \wedge 
\text{ cf}(\delta) < \lambda\}$, so $\otp(\Omega) 
\in \{\lambda^+,\lambda^+ +1\}$. 

We define $\bar e = \langle e_\delta:\delta \in \Omega\rangle$ as follows.
\bigskip

\noindent
\underline{Case 1}:  cf$(\delta) = \aleph_0,e_\delta$ is the
$<_*$-minimal member of $\{u \subseteq \delta:\delta = \sup(u)$ and
otp$(u)=0\}$.
\bigskip

\noindent
\underline{Case 2}:  cf$(\delta) > \aleph_0$.

Let $e_\delta = \cap\{c \ell_F(C):C$ a club of $\delta\}$.

So
\mn
\begin{enumerate}
\item[$(*)_1$]  $e_\delta$ is an unbounded subset of $\delta$ of order
type $< \lambda$.
\end{enumerate}
\mn
[Why?  If cf$(\delta) = \aleph_0$ then $e_\delta$ has order type
$\omega$ which is $< \lambda$ by clause (b) of the assumption.

If cf$(\delta) > \aleph_0$ then for some club $C$ of $\delta,e_\delta
= c \ell_F(C)$ has otp$(e_\delta) \le |c \ell_F(C)| \le
(\text{cf}(\delta)^{\aleph_0} < \lambda$.  The last inequality holds
as cf$(\delta) \le \lambda$ as $\delta < \lambda^+$, cf$(\delta) \ne
\lambda$ as $\lambda$ is singular by clause (c) of the assumption, and
lastly $((\text{cf}(\delta)^{\aleph_0}) < \lambda$ by clause (b) of
the assumption.]

This is enough for part (2).
Now we shall define a one-to-one function $f_\alpha$ from $\alpha$
into $\lambda$ by induction on $\alpha \in \Omega$ as
follows: let pr$_\lambda:\lambda \times \lambda \rightarrow \lambda$
be a pairing function so one to one (can add ``onto $\lambda$"); if we
succeed then $f_{\lambda^+}$ cannot be well defined so $\lambda^+
\notin \Omega$ hence cf$(\lambda^+) \ge \lambda$, but $\lambda$ is
singular so cf$(\lambda^+) = \lambda^+$, i.e.
$\lambda^+$ is not singular so we shall be done proving part (1).

The inductive definition is:
\mn
\begin{enumerate}
\item[$\boxplus$] 
\begin{enumerate}
\item[(a)]   \underline{if $\alpha \le \lambda$}
\then \, $f_\alpha$ is the identity
\sn
\item[(b)]  \underline{if $\alpha = \beta +1 \in
[\lambda,\lambda^+)$} \then \, for $i < \alpha$ we let $f_\alpha(i)$
be
\sn
\begin{itemize}
\item  $1 + f_\beta(i)$ if $i < \beta$
\sn
\item   $0$ if $i=\beta$
\end{itemize}
\sn
\item[(c)]  \underline{if $\alpha \in \Omega$} so
$\alpha$ is a limit ordinal, $e_\alpha \subseteq \alpha =
\sup(e_\alpha),e_\alpha$ of cardinality 
$< \lambda$ and we let $f_\alpha$ be defined by: for $i <\alpha$ we let 
$f_\alpha(i) = \text{ pr}_\lambda
(f_{\text{min}(e_\alpha \backslash (i+1))}(i),\otp(e_\alpha \cap i))$.
\end{enumerate}
\end{enumerate}
\end{PROOF}

\noindent
We later add:
\begin{claim}
\label{c17}
[$\ZFC$]  Assume $\mu > \kappa = \cf(\mu) > \aleph_0$ and $\mu = \mu^{\aleph_0}
+ 2^{2^\kappa}$.

\noindent
1) From some $X \subseteq \mu$ we can define well ordering of some set
$\cG \subseteq {}^\kappa \mu$ such that ${}^\kappa \mu = \{\sup\{f_n:n
< \omega\}:f_n \in \cG$ for $n < \omega\}$.

\noindent
2) If moreover $2^{2^\theta} \le \mu$ where $\theta = \kappa^\aleph$
\then \, from some $X \subseteq
\mu$ we can define a well ordering of ${}^\kappa \mu$.
\end{claim}

\begin{PROOF}{\ref{c17}}
1) Let $X \subseteq \mu$ code $\cP(\cP(\kappa)),{}^\omega \mu$ and
$F:{}^\omega \mu \rightarrow \omega$ which as in \ref{c2}.  Unlike the
proof of \ref{c2} we do not use the $\cd_i(i < \kappa)$ and we use the
family of $\aleph_1$-complete filters on $\kappa$, the rest should be
clear.

\noindent
2) As $\theta = \theta^{\aleph_0}$ there is a one-to-one onto function
$\cd:{}^{\omega}\theta \rightarrow \theta$ onto $\theta$, and for $i <
\omega$ let $\cd_i:\theta \rightarrow \theta$ be such that:
\mn
\begin{enumerate}
\item[$(*)_1$]  if $\cd(\eta) = \zeta$, \then \, $\cd_0(\zeta) = \ell
  g(\eta)$ and $\cd_{1+i}(\zeta) = \eta(i)$ for $i < \ell g(\eta)$.
\end{enumerate}
\mn
Let $D$ be $\{A \subseteq \theta$: for some $u \in [\theta]^{\le
  \aleph_0}$ we have $A \supseteq \{\varp < \theta:u \subseteq
\{\cd_i(\varp):i < \omega\}\}$, so
\mn
\begin{enumerate}
\item[$(*)_2$]  $D$ is an $\aleph_1$-complete filter on $\theta$.
\end{enumerate}
\mn
[Why?  Should be clear.]
\mn
\begin{enumerate}
\item[$(*)_3$]  for $f \in {}^\theta \mu$ let $g,g_f$ be the unique
  function $g$ with doman $\theta$ such that:
\sn
\begin{itemize}
\item  if $\varp < \kappa$ and $i < \cd_0(\varp)$, \then \,
$\cd_{1+i}(\varp) < \theta \Rightarrow 
\cd_{1+i}(g(\varp)) = f(\cd_{1+i}(\varp))$ and 
$\cd_0(g(\varp)) = \cd_0(\varp)$ and $f(\zeta) = 0$ otherwise
\end{itemize}
\end{enumerate}
\mn
[Why $g_f$ exists?  Just think.]
\mn
\begin{enumerate}
\item[$(*)_4$]  if $f \in {}^\theta \mu,\alpha = \rk_D(g_f)$ and
  $\gy = \gy_{g_f}$ as in the proof of \ref{c2} for $g_f$, \then \,:
\sn
\begin{enumerate}
\item[(a)]  from $g_f \rest Z_{\gy}$ we can define $f$ (using some $Y
  \subseteq \kappa$ as a parameter)
\sn
\item[(b)]   $\Rang(f) \subseteq \{\cd_{1+i}(g_f(\varp)):\varp \in
  Z_{\gy}$ and $i < \cd_0(g_f(\varp))\}$.
\end{enumerate}
\end{enumerate}
\mn
[Why?  Clause (a) follows clause (b).  Clause (b) holds as for every
$\xi < \kappa$, the set $\{\varp < \theta:\xi \in \{\cd_{1+i}(\varp):i
< \cd_0(\varp)\}\} \in D$.]

We continue as in the proof of \ref{c2}.
\end{PROOF}

\begin{conclusion}
\label{c20}
[$\DC$] Assume $[\lambda]^{\aleph_0}$ is well ordered for every $\lambda$.

\noindent
1) If $2^{2^\kappa}$ is well ordered \then \, for every
$\lambda,[\lambda]^\kappa$ is well ordered.

\noindent
2) For any set $Y$, there is a derived set $Y_*$ so called
{\rm Fil}$^4_{\aleph_1}(Y)$ of power near ${\cP}({\cP}(Y))$ such that
$\Vdash_{\text{\rm Levy}(\aleph_0,Y)}$ ``for every 
$\lambda,{}^Y \lambda$ is well ordered". 
\end{conclusion}

\begin{PROOF}{\ref{c20}}
1) By \ref{c2}.

\noindent
2) Follows easily.
\end{PROOF}
\newpage

\bibliographystyle{amsalpha}
\bibliography{lista,listb,listx,listf,liste,listz}

\end{document}